\documentclass[a4,semcolor,psfig,english,10pt,epsf,portrait]{article}
\usepackage{amsmath,amsfonts,latexsym,amscd,amssymb,theorem}
\usepackage{moreverb,rotating,graphics,amsfonts,latexsym,amscd}
\usepackage[ansinew]{inputenc}
\usepackage[T1]{fontenc}
\usepackage{epsfig}
\usepackage{amsmath}
\usepackage{amssymb}
\input epsf
\usepackage{psfrag}
\def\R{\hbox{\bf R}}

\def\a{\alpha}

\def\k{\kappa}

\def\eps{\varepsilon}

\def\th{\theta}
\def\r{\rho}
\def\t{\tau}

\def\<{\langle}
\def\>{\rangle}
\newcommand{\ba}{\begin{eqnarray}}
\newcommand{\ea}{\end{eqnarray}}

\newtheorem{thm}{Theorem}[section]

\newtheorem{lem}[thm]{Lemma}

\newtheorem{theorem}[thm]{Theorem}
\newtheorem{definition}[thm]{Definition}

\newtheorem{proposition}[thm]{Proposition}
\newtheorem{corollary}[thm]{Corollary}
\newtheorem{rem}[thm]{Remark}

\numberwithin{equation}{section}

\renewcommand{\R}{{\mathbb R}}

\begin{document}

\title{\bf Dynamics of dislocation densities \\in a bounded channel.
Part II: existence of weak solutions to a singular Hamilton-Jacobi/parabolic strongly
coupled system}
\author{
\normalsize\textsc{ H. Ibrahim $^{*}$, M. Jazar $^{1}$, R. Monneau
\footnote{Cermics, Paris-Est-ENPC, ParisTech, 6 et 8
avenue Blaise Pascal, Cit\'e Descartes Champs-sur-Marne, 77455
Marne-la-Vall\'ee Cedex 2, France. E-mails:
ibrahim@cermics.enpc.fr, monneau@cermics.enpc.fr.\newline \indent
${}^1$ \hskip-.1cm M. Jazar, LaMA-Liban, Lebanese University, P.O.
Box 826, Tripoli Liban, mjazar@ul.edu.lb.\newline The second
author is supported by a grant from the Lebanese University.} }}
\vspace{20pt}

\maketitle

\centerline{\footnotesize{\bf{Abstract}}} \noindent{\small{We
study a strongly coupled system consisting of a parabolic equation
and a singular Hamilton-Jacobi equation in one space dimension.
This system describes the dynamics of dislocation densities in a
material submitted to an exterior applied stress. Our system is a
natural extension of that studied in \cite{hibrahim1} where the
applied stress was set to be zero. The equations are written on a
bounded interval with Dirichlet boundary conditions and require
special attention to the boundary. We prove a result of global
existence of a solution. The method of the proof consists in
considering first a parabolic regularization of the full system,
and then passing to the limit. For this regularized system, a
result of global existence and uniqueness of a solution has been
given in \cite{IJM_PI}. We show some uniform bounds on this
solution which uses in particular an entropy estimate for the
densities.}}

\hfill\break
 \noindent{\small{\bf{AMS Classification: }}} {\small{70H20, 49L25,
      54C70, 46E30.}}\hfill\break
  \noindent{\small{\bf{Key words: }}} {\small{Hamilton-Jacobi equations,
      viscosity solutions, entropy, Orlicz spaces, parabolic equations.}}\hfill\break


\begin{section}{Introduction}\label{ch4:sec1}


\subsection{Physical motivation and setting of the problem}
In  \cite{GCZ},  Groma, Czikor and Zaiser have proposed a model
describing the dynamics of dislocation densities. Dislocations are
defects in crystals that move when a stress field is applied on
the material. These defects are one of the main explanations of
the elastoviscoplasticity behavior of metals (see
\cite{ForRad2003} and \cite{FrPiZa1993} for various models
relating dislocations and elastoviscoplastic properties of
metals). This model has been introduced in order to describe the possible
accumulation of dislocations on the boundary layer of a bounded
channel. Dislocations are distinguished by the sign of their
Burgers vector $\pm\vec{b}$ (see \cite{HL} for a description of
the Burgers vector). More precisely, let us call $\theta^{+}$ and
$\theta^{-}$, the densities of the positive and negative
dislocations respectively. For
$$x\in I:=(-1,1),$$
and $t\in(0,T)$, for some time $T>0$, the non-negative quantities $\th^{+}(x,t)$
and $\th^{-}(x,t)$ are governed by the following system (see
\cite{GCZ}):
\begin{equation}\label{ch4:sys_theta}
\left\{
\begin{aligned}
&\th^{+}_{t}=\left[\left(\frac{\th^{+}_{x}-\th^{-}_{x}}{\th^{+}+\th^{-}}
-\tau\right)\th^{+}\right]_{x}\quad&\mbox{in}&\quad I\times(0,T),\\
&\th^{-}_{t}=\left[-\left(\frac{\th^{+}_{x}-\th^{-}_{x}}{\th^{+}+\th^{-}}
-\tau\right)\th^{-}\right]_{x}\quad&\mbox{in}&\quad I\times(0,T),
\end{aligned}
\right.
\end{equation}
where $\t$ is the applied stress field which is assumed to be
constant. Here the
term $\tau_{b}=\frac{\th^{+}_{x}-\th^{-}_{x}}{\th^{+}+\th^{-}}$ is called
the back stress and can be interpreted as the contribution to the
stress of the short-range interactions between dislocations. If
$\tau=0$ and $\theta^-=0$ (resp. $\theta^+=0$) then $\theta^+$
(resp. $\theta^-$) satisfies the usual heat equation. More
generally, the back stress is proportional to the gradient of the
effective dislocation density $\theta^+-\theta^-$, with a diffusion
coefficient which is $\frac1{\theta^++\theta^-}$. In fact, system
(\ref{ch4:sys_theta}) is a model for a 2D channel with coordinates
$(x,y)$ that is invariant in the $y$-direction (see Figure \ref{IMP_fig}).
\begin{figure}[!h]
\psfrag{e1}{$e_{1}$}
\psfrag{-1}{$-1$}
\psfrag{1}{$1$}
\psfrag{e2}{$e_{2}$}
\psfrag{tau}{$\tau$}
\centering\epsfig{file=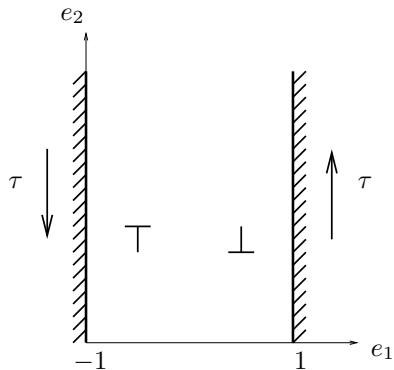, width=50mm}
   \caption{Geometry of the crystal.}\label{IMP_fig}
\end{figure}
The channel is bounded by walls that are impenetrable by
dislocations (i.e., the plastic deformation in the walls is zero).
In this case the boundary conditions are represented by the zero
flux condition, i.e.
\begin{equation}\label{ch4:cor:flux_cond}
\frac{\th^{+}_{x}-\th^{-}_{x}}{\th^{+}+\th^{-}}-\t=0,\quad
\mbox{at} \quad x=\pm1.
\end{equation}
For related literature, let us mention the work of Groma-Balogh
\cite{GB}, where the back stress was neglected. For the model
described in \cite{GB}, we refer the reader to
\cite{elhajj,elhfor} for a one-dimensional mathematical and
numerical study, and to \cite{CEMR} for a two-dimensional
existence result. The special case $\t=0$ for system
(\ref{ch4:sys_theta}) has been studied in \cite{hibrahim1}, where a
result of existence and uniqueness has been proved. In the present
paper we study the case for general constant $\t$.
\subsection{Setting of the problem}
We consider an integrated form of
(\ref{ch4:sys_theta}) and we let
$$\r^{\pm}_{x}=\th^{\pm},\quad
\r=\r^{+}-\r^{-}\quad\mbox{and}\quad\k=\r^{+}+\r^{-},$$
to obtain,  for special values of the constants of integration, the
following system in terms of $\r$ and $\k$:
\begin{equation}\label{ch4:sys_rho_kappa}
\left\{
\begin{aligned}
&\k_{t}\k_{x}=\r_{t}\r_{x}\quad&\mbox{on}&\quad I\times (0,T)\\
&\r_{t}=\r_{xx}-\tau\k_{x}\quad&\mbox{on}&\quad I\times (0,T),
\end{aligned}
\right.
\end{equation}
with the initial conditions:
\begin{equation}\label{ch4:ic}
\k(x,0)=\k^{0}(x)\quad\mbox{and}\quad\r(x,0)=\r^{0}(x).
\end{equation}
To formulate heuristically the boundary conditions at the walls
located at $x=\pm1$, we first suppose that $\k_{x}\neq 0$ at $x=\pm1$.
We recall that the dislocation fluxes at the walls must be zero, which
require (\ref{ch4:cor:flux_cond}). Rewriting system
(\ref{ch4:sys_rho_kappa}) in terms of $\rho$, $\kappa$ and
$\Phi=\frac{\th^{+}_{x}-\th^{-}_{x}}{\th^{+}+\th^{-}}-\t$, we get
\begin{equation}\label{ch4:200308b}
\left\{
\begin{aligned}
&\kappa_{t}=\rho_{x}\Phi,\\
&\rho_{t}=\kappa_x\Phi.
\end{aligned}
\right.
\end{equation}
From (\ref{ch4:cor:flux_cond}) and (\ref{ch4:200308b}), we deduce
that
\begin{equation}\label{ch4:cor:bound:r}
\r_{t}(x,.)=\k_{t}(x,.)=0 \quad \mbox{for}\quad x=\pm1.
\end{equation}
In this case, we consider the following boundary conditions:
\begin{equation}\label{ch4:bc}
\left\{
\begin{aligned}
&\k(x,.)=x\quad &\mbox{for}&\quad x=\pm 1\\
&\r(x,.)=0\quad &\mbox{for}&\quad x=\pm 1,
\end{aligned}
\right.
\end{equation}
where we have taken the zero normalization for $\r$ on the boundary of
the interval.\\

The non-negativity of $\th^{\pm}\geq 0$ reduces in terms of $\r$
and $\k$ to the following condition:
\begin{equation}\label{ch4:condition}
\k_{x}\geq |\r_{x}|,
\end{equation}
and hence a natural assumption to be considered concerning the
initial conditions $\r^{0}$ and $\k^{0}$ is to satisfy
\begin{equation}\label{ch4:cond_init}
\k^{0}_{x}\geq |\r^{0}_{x}|\quad\mbox{on}\quad I.
\end{equation}
As indicated above, problem (\ref{ch4:sys_rho_kappa}),
(\ref{ch4:ic}) and (\ref{ch4:bc}), in the case $\tau=0$, has been
studied in \cite{hibrahim1}, where a result of existence and
uniqueness is given using the viscosity/entropy solution
framework. Let us just mention that in this situation, system
(\ref{ch4:sys_rho_kappa}) becomes decoupled and easier to be
handled.


\subsection{Statement of the main result}
Remark that the first equation of system (\ref{ch4:sys_rho_kappa})
can  be formally rewritten $\kappa_t=\rho_t\rho_x/\kappa_x$ which
shows the singularity as $\kappa_x$ goes to zero. Nevertheless,
for this system we have the following result.

\begin{theorem}\textit{\textbf{(Global existence of a
solution).}}\label{ch4:mainresult}
Let $\r^{0},\k^{0}\in C^{\infty}(\bar{I})$ satisfying
(\ref{ch4:cond_init}),
\begin{equation}\label{EQNOmayma3}
\k^{0}(\pm 1)=\pm 1,\quad \r^{0}(\pm 1)=0,
\end{equation}
and the additional conditions:
\begin{equation}\label{ch4:lamo123}
D^{s}_{x}\r^{0}(\pm 1)=D^{s}_{x}\k^{0}(\pm 1)=0,\quad s=1,2.
\end{equation}
Then there exists $(\r,\k)$ such that for every
$T>0$:
$$(\r,\k)\in (C(\bar{I}\times[0,T]))^{2}\quad\mbox{and}\quad \r\in C^{1}(I\times (0,T)),$$
is a solution of (\ref{ch4:sys_rho_kappa}), (\ref{ch4:ic}) and (\ref{ch4:bc}). Moreover,
this solution satisfies (\ref{ch4:condition}) in the
distributional sense, i.e.
\begin{equation}\label{ch4:7}
\k_{x}\geq |\r_{x}|\quad\mbox{in}\quad \mathcal{D}{'}(I\times(0,T)).
\end{equation}
However, the solution has to be interpreted in the following sense:
\begin{enumerate}
  \item $\k$ is a viscosity solution of $\k_{t}\k_{x}=\r_{t}\r_{x}$ in
    $I_{T}:=I\times(0,T)$,
  \item $\r$ is a distributional solution of $\r_{t}=\r_{xx}-\tau\k_{x}$ in
  $I_{T}$,
  \item the initial and boundary conditions are satisfied
  pointwisely.
\end{enumerate}
\end{theorem}
\begin{rem}
The $C^{\infty}$ regularity of $\r^{0}$ and $\k^{0}$, together with
(\ref{ch4:lamo123}) seems to be essentially technical and are
related to an existence result for a regularized system (see Section 3,
Theorem \ref{ch4:chap3mainresult}), that we use to prove Theorem
\ref{ch4:mainresult}.
\end{rem}


\subsection{Organization of the paper}
This paper is organized as follows: in Section \ref{ch4:sec2}, we
present the strategy of the proof. In Section \ref{ch4:sec3}, we
present the tools needed throughout this work. This includes some
miscellaneous results for parabolic equations, a brief recall to
the definition and the stability  result of viscosity solutions,
and a brief recall to Orlicz spaces.  In Section
\ref{ch4:sec3bis}, we show how to choose the regularized solution.
An entropy inequality used to determine some uniform bounds on the
regularized solution is presented in Section \ref{ch4:sec4}.
Further uniform bounds and convergence arguments are done in
Section \ref{ch4:sec5}. Section \ref{ch4:sec6} is devoted to the
proof of our main result: Theorem \ref{ch4:mainresult}. In Section
\ref{ch4:Simulation}, some numerical simulations related to our
physical model are presented. Finally, Section \ref{ch4:secA} is
an appendix where we show the proofs of some technical results.
\end{section}


\begin{section}{Strategy of the proof}\label{ch4:sec2}
The main difficulty we have to face is to work with the equation
\begin{equation}\label{ch4:cor:ktkx}
\k_{t}\k_{x}=\r_{t}\r_{x}.
\end{equation}
Since $\r$ solves itself a parabolic equation (see the second equation
of (\ref{ch4:sys_rho_kappa})), we expect enough regularity on $\r$
(indeed $\r$ is $C^{1}$), and then we need a framework where the
equation involving $\k$ is stable under approximation. This property is
naturally satisfied in the framework of viscosity solutions (see
for instance \cite{barbook} and the references
therein). Then, assuming $\k_{x}\geq 0$, we interpret $\k$ as the
viscosity solution of (\ref{ch4:cor:ktkx}). Assuming
(\ref{ch4:cond_init}), we will indeed show that
$$M:=\k_{x}-|\r_{x}|\geq 0.$$
This is formally true because $M$ satisfies:
$$M_{t}=bM_{x}+cM,$$
with
$$b=\t\,\mbox{sgn}(\r_{x})-\frac{\r_{x}\r_{xx}}{\k_{x}^{2}},\quad
c=\frac{\r_{xx}^{2}}{\k_{x}^{2}}-\frac{\r_{xxx}\,\mbox{sgn}(\r_{x})}{\k_{x}},$$
where for suitable boundary conditions, we can (again formally)
see that $M\geq 0$. In order to justify the computations on $M$,
we modify the system and we consider the following parabolic
regularization for $\eps>0$ small enough:
\begin{equation}\label{ch4:app_model}
\left\{
\begin{aligned}
&\k^{\eps}_{t}=\eps\k^{\eps}_{xx}+\frac{\r^{\eps}_{x}\r^{\eps}_{xx}}{\k^{\eps}_{x}}
-\tau\r^{\eps}_{x}\quad&\mbox{in}&\quad
I\times(0,\infty)\\
&\r^{\eps}_{t}=(1+\eps)\r^{\eps}_{xx}-\tau\k^{\eps}_{x}\quad&\mbox{in}&\quad
I\times(0,\infty),
\end{aligned}
\right.
\end{equation}
 with the initial conditions:
\begin{equation}\label{EQNOmayma1}
\k^{\eps}(x,0)=\k^{0,\eps}(x),\quad \r^{\eps}(x,0)=\r^{0,\eps}(x),
\end{equation}
where $\k^{0,\eps}$ and $\r^{0,\eps}$ are some regularizations of
$\kappa^0$ and $\rho^0$ respectively,  and the same boundary conditions:
\begin{equation}\label{EQNOmayma2}
\left\{
\begin{aligned}
&\k^{\eps}(x,.)=x\quad &\mbox{for}& \quad x=\pm 1\\
&\r^{\eps}(x,.)=0\quad &\mbox{for}& \quad x=\pm 1.
\end{aligned}
\right.
\end{equation}
The system (\ref{ch4:app_model}) formally reduces to
(\ref{ch4:sys_rho_kappa}) for $\eps=0$. In
fact, system (\ref{ch4:app_model}), (\ref{EQNOmayma1}) and
(\ref{EQNOmayma2}) has (under some conditions
on the initial and boundary data) a unique
smooth global solution (see \cite[Theorem 1.1]{IJM_PI}) for $\a\in (0,1)$:
$$
(\r^{\eps},\k^{\eps})\in C^{3+\a,\frac{3+\a}{2}}(\bar{I}\times
[0,\infty))\cap  C^{\infty}(\bar{I}\times (0,\infty)).
$$
This result will be recalled in the forthcoming section (see Section 3, Theorem
\ref{ch4:chap3mainresult}). The next step is to
find some uniform bounds (independent of
$\eps$) on this solution; this is done in particular via:\\

\noindent $(1)$ an entropy inequality shown to be valid for our
regularized model (\ref{ch4:app_model});\\

\noindent $(2)$ a bound on $\k^{\eps}_{t}-\eps \k^{\eps}_{xx}$ uniformly
in $\eps$.\\

\noindent In fact, (1) guarantees the global uniform-in-time control of the
modulus of continuity in space of our regularized solution, while (2)
guarantees the local uniform-in-space control of the modulus of
continuity in time. The entropy inequality can be easily
understood. For instance, for $\eps=0$ and $\t=0$, we can formally check
that the entropy of the dislocation densities
$$\th^{\pm}=\frac{\k_{x}\pm
\r_{x}}{2},$$
defined by:
$$S(t)=\int_{I}\sum_{\pm}\th^{\pm}(.,t)\log(\th^{\pm}(.,t)),$$
satisfies:
$$\frac{dS(t)}{dt}=-\int_{I}\frac{(\th^{+}_{x}-\th^{-}_{x})^{2}}{\th^{+}+\th^{-}}\leq
0,$$
therefore we get $S(t)\leq S(0)$ which controls the entropy uniformly in
time. Finally, we need to pass to the limit $\eps\rightarrow 0$ after
multiplying the first equation of
(\ref{ch4:app_model}) by $\k^{\eps}_{x}$. Having enough control on the
regularized solutions, we can find a solution of the limit equation
using in particular the stability of viscosity solutions of
Hamilton-Jacobi equations. However, the passage to the limit in the second equation
of (\ref{ch4:app_model}) is done in the distributional sense.
\end{section}


\begin{section}{Tools: miscellaneous results on parabolic equations, viscosity
    solution, and Orlicz spaces}\label{ch4:sec3}


\subsection{Miscellaneous results on parabolic equations}

We first fix some notations. Denote
$$I_{T}:=I\times (0,T),\quad \overline{I_{T}}:=\bar{I}\times [0,T],
\quad \partial^{p}I_{T}:=I\cup(\partial I\times [0,T]),\quad
\mbox{and}\quad \|\cdot\|_{L^{p}(\Omega)}=\|\cdot\|_{p,\Omega}.$$
Define the parabolic Sobolev space $W^{2,1}_{p}(I_{T})$ , $1<p<\infty$ by:
$$W^{2,1}_{p}(I_{T}):= \left\{u\in L^{p}(I_{T});\;
  \left(u_{t},u_{x},u_{xx}\right)\in \left(L^{p}(I_{T})\right)^{3}\right\}.$$
We start with a result of global existence and uniqueness of smooth
solutions of the regularized system (\ref{ch4:app_model}), with the initial
and boundary conditions (\ref{EQNOmayma1}) and (\ref{EQNOmayma2}).
\begin{theorem}\textit{\textbf{(Global existence for the regularized system,
    \cite[Theorem 1.1]{IJM_PI}).}}\label{ch4:chap3mainresult}
Let $0<\a<1$ and $0<\eps<1$. Let $\r^{0,\eps}$, $\k^{0,\eps}$ satisfying:
\begin{equation}\label{ch4:MT1}
\r^{0,\eps},\k^{0,\eps}\in C^{\infty}(\bar{I}),\quad
\r^{0,\eps}(\pm1)=0,\quad \mbox{and}\quad \k^{0,\eps}(\pm 1)=\pm 1,
\end{equation}
\begin{equation}\label{ch4:MT2}
\left\{
\begin{aligned}
&(1+\eps)\r^{0,\eps}_{xx}=\t\k^{0,\eps}_{x}\quad &\mbox{on}&\quad \partial I\\
&(1+\eps)\k^{0,\eps}_{xx}=\t\r^{0,\eps}_{x}\quad &\mbox{on}&\quad \partial I,
\end{aligned}
\right.
\end{equation}
and
\begin{equation}\label{ch4:MT3}
\k^{0,\eps}_{x}>|\r_{x}^{0,\eps}|\quad \mbox{on}\quad \bar{I}.
\end{equation}
Then there exists a unique global solution
\begin{equation}\label{ch4:MT3nos}
(\r^{\eps},\k^{\eps})\in C^{3+\a,\frac{3+\a}{2}}(\bar{I}\times
[0,\infty))\cap  C^{\infty}(\bar{I}\times (0,\infty)),
\end{equation}
of the system  (\ref{ch4:app_model}), (\ref{EQNOmayma1}) and
(\ref{EQNOmayma2}). Moreover,
this solution satisfies :
\begin{equation}\label{ch4:MT5}
\k^{\eps}_{x}>|\r^{\eps}_{x}| \quad \mbox{on} \quad  \bar{I}\times [0,\infty).
\end{equation}
\end{theorem}
\begin{rem}
Conditions (\ref{ch4:MT2}) are natural here. Indeed, the regularity
(\ref{ch4:MT3nos}) of the solution of equation (\ref{ch4:app_model}) with
boundary conditions (\ref{EQNOmayma1}) and (\ref{EQNOmayma2}) imply in particular
condition (\ref{ch4:MT2}).
\end{rem}
\begin{rem}\textit{\textbf{(Uniform $L^{\infty}$ bound on $\r^{\eps}$ and
      $\k^{\eps}$).}}\label{ch4:fshA}
We remark, from the boundary conditions (\ref{EQNOmayma2}) and from the
inequality (\ref{ch4:MT5}), that:
\begin{equation}\label{ch4:fshA1}
\|\r^{\eps}\|_{L^{\infty}(\bar{I}\times
[0,\infty))}\leq 1\quad \mbox{and} \quad \|\k^{\eps}\|_{L^{\infty}(\bar{I}\times
[0,\infty))}\leq 1.
\end{equation}
\end{rem}
We now present two technical lemmas that will be used in the proof of
Theorem \ref{ch4:mainresult}. The proofs of these lemmas will be given
in the Appendix.
\begin{lem}\textit{\textbf{(Control of the modulus of continuity in
time uniformly in $\eps$).}}\label{ch4:mod_con_t}
Let $p>3$, and $u^{\eps}\in W^{2,1}_{p}(I_{T})$. Suppose furthermore that the sequences
$$
(u^{\eps})_{\eps}\quad\mbox{and}\quad(f^{\eps})_{\eps}=(u^{\eps}_{t}-\eps
u^{\eps}_{xx})_{\eps},
$$
are locally bounded in $I_{T}$ uniformly for $\eps\in (0,1)$. Then
for every $V\subset\!\subset I_{T}$, there exist two constants
$c>0$, $\eps_{0}>0$ depending on $V$, and $0<\beta<1$ such that for
all $0<\eps<\eps_{0}$:
$$
\frac{|u^{\eps}(x,t+h)-u^{\eps}(x,t)|}{h^{\beta}}\leq c,\quad
\forall (x,t), (x,t+h)\in V.
$$
\end{lem}
\begin{lem}\textit{\textbf{(An interior estimate for the heat
      equation).}}\label{ch4:vnlemma}
Let $a\in C^{\infty}(I_{T})\cap L^{1}(I_{T})$ satisfying:
$$
a_{t}=a_{xx} \quad \mbox{on} \quad I_{T},
$$
then for any $V \subset\!\subset I_{T}$, an open set, we have:
$$
\|a\|_{p,V}\leq c \|a\|_{1,I_{T}},\quad \forall \,1<p<\infty,
$$
where $c=c(p,V)>0$ is a positive constant.
\end{lem}

\subsection{Viscosity solution: definition and stability result}
Let $\Omega\subset\mathbb{R}^{n}$ be an open domain, and consider
the following Hamilton-Jacobi equation:
\begin{equation}\label{ch4:one}
H(x,u(x),Du(x),D^{2}u(x))=0, \quad \forall x\in \Omega,
\end{equation}
where $H:\Omega\times\mathbb{R}\times\mathbb{R}^{n}\times
M^{n\times n}_{{sym}}\rightarrow\mathbb{R}$ is a continuous
mapping.
\begin{definition}\textit{\textbf{(Viscosity solution of Hamilton-Jacobi equations).}}
A continuous function $u\,:\,\Omega\mapsto\mathbb{R}$ is a viscosity
sub-solution of (\ref{ch4:one}) if for any $\phi\in C^{2}(\Omega;\R)$
and any local maximum $x_0\in\Omega$ of $u-\phi$, one has
$$H(x_{0},u(x_{0}),D\phi(x_{0}),D^{2}\phi(x_{0}))\leq 0.$$
Similarly, $u$ is a viscosity super-solution of (\ref{ch4:one}), if at any local minimum
point $x_0\in\Omega$ of $u-\phi$,
one has
$$ H(x_{0},u(x_{0}),D\phi(x_{0}),D^{2}\phi(x_{0}))\geq 0.$$
Finally, if $u$ is both a
viscosity sub-solution and a viscosity super-solution, then $u$ is
called a viscosity solution.
\end{definition}
To get a "non-empty" and useful definition, it is usually assumed
that $H$ is elliptic (see \cite{barbook}). This notion of
ellipticity will be indirectly used in Section \ref{ch4:sec6}. In
fact, this definition is used for interpreting solutions of the
first equation of (\ref{ch4:sys_rho_kappa}) in the viscosity
sense. This will be shown in Section \ref{ch4:sec4}. To be more
precise, in the case where $\Omega=I_{T}$, we say that $u$ is a
viscosity solution of the Dirichlet problem (\ref{ch4:one}) with
$u=\zeta\in C(\partial^{p}I_{T})$ if:
\begin{description}
  \item[] (1) $u\in C(\overline{I_{T}})$,
  \item[] (2) $u$ is a viscosity solution of (\ref{ch4:one}) in
  $I_{T}$,
  \item[] (3) $u=\zeta$ on $\partial^{p}I_{T}$.
\end{description}
For a better understanding of the viscosity interpretation of
boundary conditions of Hamilton-Jacobi equations, we refer the
reader to \cite[Section 4.2]{barbook}. We now state a stability
result for viscosity solutions of Hamilton-Jacobi equations.
\begin{theorem}\textit{\textbf{(Stability of viscosity solutions,  \cite[Lemma
    2.3]{barbook}).}}\label{ch4:stability}
Suppose that, for $\eps>0$, $u^{\eps}\in C(\Omega)$ is a viscosity
sub-solution (resp. super-solution) of the equation
$$
H^{\eps}(x,u^{\eps},Du^{\eps},D^{2}u^{\eps})=0\quad \mbox{in}\quad
\Omega,
$$
where $(H^{\eps})_{\eps}$ is a sequence of continuous functions. If
$u^{\eps}\rightarrow u$ locally uniformly in $\Omega$ and if
$H^{\eps}\rightarrow H$ locally uniformly in
$\Omega\times\R\times\R^{n}\times M^{n\times n}_{sym}$, then $u$ is
a viscosity sub-solution (resp. super-solution) of the equation:
$$
H(x,u,Du,D^{2}u)=0\quad\mbox{in}\quad \Omega.
$$
\end{theorem}


\subsection {Orlicz spaces: definition and properties}
We recall the definition of an Orlicz space and some of its
properties (for details see \cite{adams}). A real valued function
$\Psi:[0,\infty)\rightarrow \R$ is called a Young function if
$$\Psi(t)=\int_{0}^{t}\psi(s)ds,$$
where $\psi:[0,\infty)\rightarrow [0,\infty)$ satisfying:
\begin{itemize}
  \item $\psi(0)=0$, $\psi>0$ on $(0,\infty)$, $\psi(t)\rightarrow
  \infty$ as $t\rightarrow \infty$;
  \item $\psi$ is non-decreasing and right continuous at any point
  $s\geq 0$.
\end{itemize}
Let $\Psi$ be a Young function. The Orlicz class $K_{\Psi}(I)$ is
the set of equivalence classes of real-valued measurable functions
$u$ on $I$ satisfying
$$\int_{I}\Psi(|u(x)|)dx< +\infty.$$
\begin{definition}\textit{\textbf{(Orlicz spaces).}}\label{ch4:orlicz spaces}
The Orlicz space $L_{\Psi}(I)$ is the linear span of $K_{\Psi}(I)$
supplemented with the Luxemburg norm
\begin{equation}\label{ch4:norm_O}
\|u\|_{L_{\Psi}(I)}=\inf\left\{k>0;\;\int_{I}\Psi\left(\frac{|u(x)|}{k}\right)\leq
1\right\},
\end{equation}
and with this norm, the Orlicz space is a Banach space.
\end{definition}
The function
$$\Phi(t)=\int_{0}^{t}\phi(s)ds,\quad \phi(s)=\sup_{\psi(t)\leq s}t,$$
is called the complementary Young function of $\Psi$. An example of
such pair of complementary Young functions is the following:
\begin{equation}\label{ch4:Young_comp}
\Psi(s)=(1+s)\log(1+s)-s\quad\mbox{and}\quad\Phi(s)=e^{s}-s-1.
\end{equation}
We now state a lemma giving two useful properties of Orlicz spaces that
will be used in the proof of Lemma \ref{ch4:control_x}.
\begin{lem}\textit{\textbf{(Norm control and Hölder inequality,
      \cite{KJF}).}}\label{ch4:cor:lem:orli}
If $u\in L_{\Psi}(I)$ for some Young function $\Psi$, then we have:
\begin{equation}\label{ch4:prop1}
\|u\|_{L_{\Psi}(I)}\leq 1+ \int_{I}\Psi(|u(x)|)dx.
\end{equation}
Moreover, if $v\in L_{\Phi}(I)$, $\Phi$ being the complementary Young
function of $\Psi$, then we have the following Hölder inequality:
\begin{equation}\label{ch4:prop2}
\left|\int_{I}uv dx\right|\leq
2\|u\|_{L_{\Psi}(I)}\|v\|_{L_{\Phi}(I)}.
\end{equation}
\end{lem}
\end{section}


\begin{section}{The regularized problem}\label{ch4:sec3bis}
As we have already mentioned, we will use a parabolic regularization of
(\ref{ch4:sys_rho_kappa}), and a result of global existence of this
regularized system from \cite{IJM_PI} (see Theorem
\ref{ch4:chap3mainresult}). In order to use this result, we need to give
a special attention to the conditions on the initial data of the
approximated system $\r^{0,\eps}$ and $\k^{0,\eps}$ (see
(\ref{ch4:MT1}), (\ref{ch4:MT2}) and (\ref{ch4:MT3})). This section aims
to show how to choose the suitable initial data $\r^{0,\eps}$ and
$\k^{0,\eps}$ in order to benefit Theorem \ref{ch4:chap3mainresult}.
Let $\r^{0}$ and $\k^{0}$ be the functions given in Theorem \ref{ch4:mainresult}.
Set
\begin{equation}\label{ch4:cond_init_ap2}
\r^{0,\eps}(x)=\frac{\r^{0}(x)+\eps\t\phi(x)}{(1+\eps)^{2}}\quad
\mbox{and}\quad \k^{0,\eps}(x)=\frac{\k^{0}(x)+\eps x}{1+\eps},
\end{equation}
with the function $\phi$ defined by:
\begin{equation}\label{ch4:functionphi}
\phi(x)=\left\{\begin{array}{ll}\frac{1}{4\t^{2}}[1-\cos\t(x^{2}-1)]&\mbox{
      if } \,\tau\ne0\\ 0&\mbox{ if }\,\tau=0.\end{array}\right.
\end{equation}
The function $\phi$ enjoys some properties that are shown in the
following lemma.
\begin{lem}\textit{\textbf{(Properties of $\phi$)}}\label{ch4:RemarK}\\
The function $\phi$ given by (\ref{ch4:functionphi}) satisfies the following
properties:\\

(P1) $\phi,\phi^{'}|_{\partial I}=0$;\\

(P2) $\phi^{''}\big|_{\partial I}=1$;\\

(P3) $\displaystyle |\phi^{'}(x)|< {1}/{|\t|}\quad$ \rm{for} $\quad x\in \bar{I}$.
\end{lem}
\noindent {\bf Proof.} Straightforward computations. $\hfill{\Box}$\\

Form the above lemma, and from the construction of $\r^{0,\eps}$ and
$\k^{0,\eps}$ (see (\ref{ch4:cond_init_ap2})) together with the properties enjoyed by $\r^{0}$
and $\k^{0}$ (see (\ref{EQNOmayma3}) and (\ref{ch4:lamo123})), we write down some properties of
$\r^{0,\eps}$ and $\k^{0,\eps}$.
\begin{lem}\textit{\textbf{(Properties of $\r^{0,\eps}$ and $\k^{0,\eps}$)}}\label{ch4:RemarK1}\\
The functions  $\r^{0,\eps}$ and $\k^{0,\eps}$ given by
(\ref{ch4:cond_init_ap2})  , satisfy the following
properties:\\

(P4) $\r^{0,\eps}(\pm1)=\pm 1,\quad \mbox{and}\quad \k^{0,\eps}(\pm 1)=0$; \\

(P5) $\displaystyle (1+\eps)\k^{0,\eps}_{xx}\big|_{\partial
  I}=\t\r^{0,\eps}_{x}\big|_{\partial I}\quad $and $\quad (1+\eps)\r^{0,\eps}_{xx}\big|_{\partial
  I}=\t\k^{0,\eps}_{x}\big|_{\partial I}$;\\

(P6) $\displaystyle \k^{0,\eps}_{x}\geq
|\r^{0,\eps}_{x}|+\frac{\eps (1-|\tau||\phi^{'}|)}{1+\eps}>|\r^{0,\eps}_{x}|$.
\end{lem}
\noindent {\bf Proof.} Straightforward computations.
$\hfill{\Box}$

\begin{rem}\textit{\textbf{(The regularized solution
      $(\r^{\eps},\k^{\eps})$).}}\label{ch4:olia}
Properties (P4)-(P5)-(P6) of Lemma \ref{ch4:RemarK1} implies condition
(\ref{ch4:MT1})-(\ref{ch4:MT2})-(\ref{ch4:MT3}) of Theorem
\ref{ch4:chap3mainresult}. In this case, call
\begin{equation}\label{ch4:nnifix}
(\r^{\eps},\k^{\eps}),
\end{equation}
the solution of
(\ref{ch4:app_model}), (\ref{EQNOmayma1}) and
(\ref{EQNOmayma2}), given in Theorem \ref{ch4:chap3mainresult}, with the initial
conditions
$$
\r(x,0)=\r^{0,\eps}\quad \mbox{and}\quad \k(x,0)=\k^{0,\eps},
$$
that are given by (\ref{ch4:cond_init_ap2}).
\end{rem}
\end{section}


\begin{section}{Entropy inequality}\label{ch4:sec4}
\begin{proposition}\textit{\textbf{(Entropy inequality).}}\label{ch4:ent_ineq}
Let $(\r^{\eps},\k^{\eps})$ be the regular solution given by (\ref{ch4:nnifix}).
Define $\theta^{\pm,\eps}$ by:
\begin{equation}\label{ch4:thetapm}
\theta^{\pm,\eps}=\frac{\k^{\eps}_{x}\pm \r^{\eps}_{x}}{2},
\end{equation}
then the quantity $S(t)$ given by:
\begin{equation}\label{ch4:S(t)}
S(t)=\int_{I}\sum_{\pm}\theta^{\pm,\eps}(x,t)\log \theta^{\pm,\eps}(x,t) dx,
\end{equation}
satisfies for every $t\geq 0$:
\begin{equation}\label{ch4:eqn_sat_S(t)}
S(t)\leq S(0)+ \frac{\tau^{2}t}{2} .
\end{equation}
\end{proposition}
\noindent {\bf Proof.} From (\ref{ch4:MT5}), we know that
$\k^{\eps}_{x}>|\r^{\eps}_{x}|$, hence $\th^{\pm,\eps}>0$,
and the term $\log(\th^{\pm,\eps})$ is well defined. Also from the
regularity (\ref{ch4:MT3nos}) of the solution $(\r^{\eps},\k^{\eps})$,
we know that $\th^{\pm,\eps}(.,t)\in C(\bar{I})$ for all  $t\geq 0$,
hence the term $S(t)$ is well defined. We derive system
(\ref{ch4:app_model}) with respect to $x$, and we write it in terms of
$\theta^{\pm,\eps}$, we get:
\begin{equation}\label{ch4:conservation}
\left\{
\begin{aligned}
&\theta^{+,\eps}_{t}=\left[\left(\frac{(\theta^{+,\eps}-\theta^{-,\eps})_{x}}{\theta^{+,\eps}+\theta^{-,\eps}}-\tau
\right)\theta^{+,\eps}+\eps\theta^{+,\eps}_{x}\right]_{x}\\
&\theta^{-,\eps}_{t}=\left[-\left(\frac{(\theta^{+,\eps}-\theta^{-,\eps})_{x}}{\theta^{+,\eps}+\theta^{-,\eps}}-\tau
\right)\theta^{-,\eps}+\eps\theta^{-,\eps}_{x}\right]_{x}.
\end{aligned}
\right.
\end{equation}
We first remark that:
$$\pm\left(\frac{(\theta^{+,\eps}-\theta^{-,\eps})_{x}}{\theta^{+,\eps}+\theta^{-,\eps}}-\tau
\right)\theta^{\pm,\eps}+\eps\theta^{\pm,\eps}_{x}=\frac{\k^\eps_{t}\pm\r^\eps_{t}}{2}.$$
Since $\k^{\eps}_{t}$ and $\r^{\eps}_{t}$ are zeros on $\partial I
\times [0,\infty)$, then
\begin{equation}\label{ch4:zero_b_ent}
\left(\frac{(\theta^{+,\eps}-\theta^{-,\eps})_{x}}{\theta^{+,\eps}+\theta^{-,\eps}}-\tau
\right)\theta^{+,\eps}+\eps\theta^{+,\eps}_{x}=-\left(\frac{(\theta^{+,\eps}-
\theta^{-,\eps})_{x}}{\theta^{+,\eps}+\theta^{-,\eps}}-\tau
\right)\theta^{-,\eps}+\eps\theta^{-,\eps}_{x}=0\;\mbox{on}\;\partial I\times [0,\infty).
\end{equation}
Using (\ref{ch4:zero_b_ent}), we compute
for $t\geq 0$:
\begin{eqnarray*}
S^{'}(t)&=&\sum_{\pm}\int_{I}\theta^{\pm,\eps}_{t}\log(\theta^{\pm,\eps})+\theta^{\pm,\eps}_{t},\\
&=&\sum_{\pm}\int_{I}\mp\left(\frac{(\theta^{+,\eps}-\theta^{-,\eps})_{x}}{\theta^{+,\eps}+\theta^{-,\eps}}-\tau
\right)\theta^{\pm,\eps}_{x}-\eps\frac{\big(\theta^{\pm,\eps}_{x}\big)^{2}}{\theta^{\pm,\eps}},\\
&=&\int_{I}-\frac{\big(\theta^{+,\eps}_{x}-\theta^{-,\eps}_{x}\big)^{2}}{\theta^{+,\eps}+\theta^{-,\eps}}
+\tau(\theta^{+,\eps}_{x}-\theta^{-,\eps}_{x})
-\eps\left(\frac{\big(\theta^{+,\eps}_{x}\big)^{2}}{\theta^{+,\eps}}+
\frac{\big(\theta^{-,\eps}_{x}\big)^{2}}{\theta^{-,\eps}}\right),
\end{eqnarray*}
where we have integrated by parts in the second line. By Young's
Inequality, we have:
$$
\left|\th^{+,\eps}_{x}-\th^{-,\eps}_{x}\right|\leq
\frac{1}{|\tau|}\frac{\big(\th^{+,\eps}_{x}-\th^{-,\eps}_{x}\big)^{2}}{\th^{+,\eps}+\th^{-,\eps}}+
\frac{|\tau|}{4}(\th^{+,\eps}+\th^{-,\eps}),
$$
and hence
$$S^{'}(t)\leq\frac{\tau^{2}}{4}\int_{I}(\th^{+,\eps}+\theta^{-,\eps}).$$
Moreover, we have from (\ref{EQNOmayma2}), that
$$\int_{I}(\th^{+,\eps}(.,t)+\th^{-,\eps}(.,t))=\int_{I}\k^{\eps}_{x}(.,t)=\k^{\eps}(1,t)-\k^{\eps}(-1,t)=2,$$
and therefore
$$S^{'}(t)\leq \frac{\tau^{2}}{2}.$$
Integrating the previous inequality from $0$ to $t$, we get
(\ref{ch4:eqn_sat_S(t)}). $\hfill{\Box}$\\

\noindent An immediate corollary of Proposition \ref{ch4:ent_ineq} is the
following:
\begin{corollary}\textit{\textbf{(Special control of  $\k^{\eps}_{x}$)}}\label{ch4:k_{x}}\\
For all $t\geq 0$, we have:
\begin{equation}\label{ch4:cont_k_x}
\int_{I}\k^{\eps}_{x}(x,t)\log(\k^{\eps}_{x}(x,t))dx\leq
S(0)+\frac{\tau^{2}t}{2}+2,
\end{equation}
where $S$ is given by (\ref{ch4:S(t)}).
\end{corollary}

\noindent The proof of Corollary \ref{ch4:k_{x}} comes from the following
inequality.

\begin{lem}\label{ch4:usef_ineq1}
For every $x,y>0$, we have:
\begin{equation}\label{ch4:usef_ineq2}
(x+y)\log(x+y)\leq x\log(x)+y\log(y)+x\log(2)+y.
\end{equation}
\end{lem}
\noindent {\bf Proof.} Direct computations. $\hfill{\Box}$\\

\noindent {\bf Proof of Corollary \ref{ch4:k_{x}}.} From (\ref{ch4:thetapm}),
it follows that
$$\k^{\eps}_{x}=\th^{+,\eps}+\th^{-,\eps}>0.$$
Then we have for $t\geq 0$:
\begin{eqnarray*}
\int_{I}\k^{\eps}_{x}\log \k^{\eps}_{x} &=&\int_{I}(\th^{+,\eps}+\th^{-,\eps})\log(\th^{+,\eps}+\th^{-,\eps})\\
&\leq&
\int_{I}\th^{+,\eps}\log(\th^{+,\eps})+\th^{-,\eps}\log(\th^{-,\eps})+\th^{+,\eps}\log2+\th^{-,\eps}\\
&\leq&
\left(\int_{I}\th^{+,\eps}\log(\th^{+,\eps})+\th^{-,\eps}\log(\th^{-,\eps})\right)+\log2+1\\
&\leq& S(t)+2.
\end{eqnarray*}
Here we have used Lemma \ref{ch4:usef_ineq1} with $x=\th^{+,\eps}$ and
$y=\th^{-,\eps}$ for the second line, and we have used for the third
line, the fact that
$$\int_{I}\th^{\pm,\eps}=\frac{1}{2}\int_{I}\k^{\eps}_{x}\pm
\r^{\eps}_{x}=\frac12[\k^{\eps}(1,.)-\k^{\eps}(-1,.)]=1.$$
Using (\ref{ch4:eqn_sat_S(t)}), the result follows. $\hfill{\Box}$

\begin{lem}\textit{\textbf{(Control of the modulus of continuity in
space)}}\label{ch4:control_x}\\
Let $u\in C^{1}(I)$, $u_{x}>0$, satisfying
$$
\int_{I}u_{x}\log(u_{x})\leq c_{1},
$$
for some positive constant $c_1$, then we have for any $x$,
$x+h\in I$ with $h>0$:
\begin{equation}\label{ch4:bol}
|u(x+h)-u(x)|\leq \frac{2(c_{1}+1+\log2)}{|\log h|}.
\end{equation}
\end{lem}

\noindent {\bf Proof.} Let $x, x+h\in I$.\\

\noindent {\bf Step 1. }\textsf{($u_{x}\in L_{\Psi}(x,x+h)$ with $\Psi$
  given in (\ref{ch4:Young_comp}))}\\

\noindent  We compute
\begin{eqnarray*}
\int_{x}^{x+h}\Psi(u_{x}) &\leq& \int_{I}(1+u_{x})\log(1+u_{x})-u_{x}\\
&\leq& \int_{I} u_{x}\log(u_{x})+\log 2\leq c_{1}+\log 2,
\end{eqnarray*}
where we have used (\ref{ch4:usef_ineq2}) in the second line
inequality. Hence from (\ref{ch4:prop1}), we get
$$\|u_{x}\|_{L_{\Psi}(x,x+h)}\leq c_{1}+1+\log 2.$$

\noindent {\bf Step 2. }\textsf{(Estimating the modulus of
  continuity)}\\

\noindent It is easy to check that the function $1$ lies in $
L_{\Phi}(x,x+h)$ for $\Phi$ given by (\ref{ch4:Young_comp}), and that
$\|1\|_{L_{\Phi}(x,x+h)}\leq~-\frac{1}{\log h}$.
Therefore, by Hölder inequality (\ref{ch4:prop2}), we obtain:
\begin{eqnarray*}
|u(x+h)-u(x)|&=&\left|\int_{x}^{x+h}u_{x}.1\right|\nonumber\\
&\leq & 2\|u_{x}\|_{L_{\psi}(x,x+h)}\|1\|_{L_{\Phi}(x,x+h)}\leq
\frac{2(c_{1}+1+\log 2)}{|\log h|},
\end{eqnarray*}
and the result follows. $\hfill{\Box}$
\begin{rem}
As mentioned to us by J\'{e}r\^{o}me Droniou, it is possible to estimate
directly the quantity $|u(x+h)-u(x)|\leq \frac{A}{|\log h|}$ by
splitting the integral
$\int^{x+h}_{x}u_{x}$ on the set where $u_{x}$ is bigger and lower than
$\lambda$, and then optimizing on the parameter $\lambda$.
\end{rem}
\end{section}


\begin{section}{An interior estimate}\label{ch4:sec5}
In this section, we give an interior estimate for the term
\begin{equation}\label{ch4:qA}
A^{\eps}=\r^{\eps}_{x}-\tau\k^{\eps}.
\end{equation}
that will be used in the passage to the limit as $\eps$ goes to
zero in the regularized system. We start by deriving an equation
satisfied by $A^{\eps}$:
\begin{equation}\label{ch4:para_qA}
A^{\eps}_{t}=(1+\eps)A^{\eps}_{xx}-\frac{\tau\r^{\eps}_{x}}{\k^{\eps}_{x}}A^{\eps}_{x}.
\end{equation}

\noindent We now show an interior $L^{p}$ estimate concerning
the term $A^{\eps}$. This estimate gives a control on the local $L^{p}$ norm
of $A^{\eps}$ by its global $L^{1}$ norm over $I_{T}$, and it will be
used in the following section. More precisely, we have the
following lemma.
\begin{lem}\textit{\textbf{(Interior $L^{p}$ estimate).}}\label{ch4:3lcom:lem1}
Let $0<\eps<1$ and $1<p<\infty$. Then the quantity $A^{\eps}$ given by
(\ref{ch4:qA}) satisfies:
\begin{equation}\label{ch4:3:1}
\|A^{\eps}\|_{p,V}\leq c \left(\|A^{\eps}\|_{1,I_{T}}+1 \right),
\end{equation}
where $V$ is an open subset of $I_{T}$ such that $V\subset\!\subset
I_{T}$, and $c=c(p,V)>0$ is a constant independent of $\eps$.
\end{lem}
{\bf Proof.} Throughout the proof, the term $c=c(p,V)>0$ is a
positive constant independent of $\eps$, and it may vary from line to
line. A simple computation gives:
\begin{eqnarray}\label{ch4:aymta5las}
-\tau\frac{\r^{\eps}_{x}}{\k^{\eps}_{x}}A^{\eps}_{x}&=&
-\tau\frac{\r^{\eps}_{x}}{\k^{\eps}_{x}}(\r^{\eps}_{xx}-\tau\k^{\eps}_{x})\nonumber\\
&=&
-\tau\frac{\r^{\eps}_{x}\r^{\eps}_{xx}}{\k^{\eps}_{x}}+\tau^{2}\r^{\eps}_{x}\nonumber\\
&=&-\tau (\k^{\eps}_{t}-\eps\k^{\eps}_{xx}).
\end{eqnarray}
Define $\bar{\k}^{\eps}$ as the unique solution of
\begin{equation}\label{ch4:3l:1}
\left\{
\begin{aligned}
&\bar{\k}^{\eps}_{t}=(1+\eps)\bar{\k}^{\eps}_{xx}+\k^{\eps}\quad
&\mbox{on}&\quad I_{T},\\
&\bar{\k}^{\eps}=0 \quad&\mbox{on}& \quad \partial^{p}I_{T},
\end{aligned}
\right.
\end{equation}
where the existence and uniqueness of this equation is a direct
consequence of the $L^{p}$ theory for parabolic equations (see for
instance \cite[Theorem 9.1]{LSU}) using in particular the fact
that $\k^{\eps}\in C^{1}(\overline{I_{T}})$. Moreover, from the
regularity (\ref{ch4:MT3nos}) of $\k^{\eps}$, we can deduce that
$\bar{\k}^{\eps}\in C^{\infty}(I_{T})$. Let
$$\bar{A}^{\eps}=-\tau
(\bar{\k}^{\eps}_{t}-\eps\bar{\k}^{\eps}_{xx}),\quad \mbox{with} \quad
a^{\eps}=A^{\eps}-\bar{A}^{\eps}.$$
We calculate:
\begin{eqnarray*}
\bar{A}^{\eps}_{t}&=&-\tau[\bar{\k}^{\eps}_{tt}-\eps
\bar{\k}^{\eps}_{xxt}]\\
&=& -\tau[(1+\eps)\bar{\k}^{\eps}_{xxt}+\k^{\eps}_{t}-\eps
((1+\eps)\bar{\k}^{\eps}_{xxxx}+\k^{\eps}_{xx})]\\
&=&-\tau(1+\eps)(\bar{\k}^{\eps}_{xxt}-\eps\bar{\k}^{\eps}_{xxxx} )
-\tau (\k^{\eps}_{t}-\eps\k^{\eps}_{xx} )\\
&=&(1+\eps)\bar{A}_{xx}-\frac{\tau\r^{\eps}_{x}}{\k^{\eps}_{x}}A^{\eps}_{x},
\end{eqnarray*}
where for the first two line, we have used (\ref{ch4:3l:1}), and for the
last line, we have used (\ref{ch4:aymta5las}). In this case, we obtain:
\begin{eqnarray*}
a^{\eps}_{t}&=& A^{\eps}_{t}-\bar{A}^{\eps}_{t}\\
&=&
(1+\eps)A^{\eps}_{xx}-\frac{\tau\r^{\eps}_{x}}{\k^{\eps}_{x}}A^{\eps}_{x}-(1+\eps)\bar{A}_{xx}
+\frac{\tau\r^{\eps}_{x}}{\k^{\eps}_{x}}A^{\eps}_{x}\\
&=& (1+\eps)
(A^{\eps}_{xx}-\bar{A}^{\eps}_{xx})=(1+\eps)a^{\eps}_{xx},
\end{eqnarray*}
where for the first line, we have used the equation
(\ref{ch4:para_qA}). We apply Lemma \ref{ch4:vnlemma} to the
function $a^{\eps}$, after doing parabolic rescaling of the form
$\tilde{a}^{\eps}(x,t)=a^{\eps}\left(x,\frac{t}{1+\eps}\right)$,
we get:
$$
\|a^{\eps}\|_{p,V}\leq c(1+\eps)^{1-\frac{1}{p}} \|a^{\eps}\|_{1,I_{T}},
$$
and since $0<\eps<1$, we finally obtain
$$
\|a^{\eps}\|_{p,V}\leq c\|a^{\eps}\|_{1,I_{T}}.
$$
From the definition of $a^{\eps}$, and the last inequality, we
finally deduce that:
\begin{equation}\label{ch4m7styet}
\|A^{\eps}\|_{p,V}\leq c
(\|A^{\eps}\|_{1,I_{T}}+\|\bar{A}^{\eps}\|_{p,I_{T}}).
\end{equation}
In order to complete the proof, we need to control the term
$\|\bar{A}^{\eps}\|_{p,I_{T}}$ in (\ref{ch4m7styet}). We use the
equation (\ref{ch4:3l:1}) satisfied by $\bar{\k}^{\eps}$ to
obtain:
\begin{eqnarray}\label{ch4:RbT}
\|\bar{A}^{\eps}\|_{p,I_{T}}&=&\tau \|\bar{\k}^{\eps}_{t}-\eps
\bar{\k}^{\eps}_{xx}\|_{p,I_{T}}\nonumber\\
&=& \tau \|\bar{\k}^{\eps}_{xx}+\k^{\eps}\|_{p,I_{T}}\nonumber\\
&\leq& c(\|\bar{\k}^{\eps}_{xx}\|_{p,I_{T}}+\|\k^{\eps}\|_{p,I_{T}}).
\end{eqnarray}
The $L^{p}$ estimates for parabolic equations (see \cite[Lemma
2.7]{IJM_PI}) applied to (\ref{ch4:3l:1}) gives:
$$\|\bar{\k}^{\eps}_{xx}\|_{p,I_{T}}\leq
\frac{c}{1+\eps}\|\k^{\eps}\|_{p,I_{T}},$$
then (\ref{ch4:RbT}), together with the fact that $0\leq \k^{\eps} \leq
1$ (see Remark \ref{ch4:fshA}), implies that:
$$\|\bar{A}^{\eps}\|_{p,I_{T}}\leq c \|\k^{\eps}\|_{p,I_{T}}\leq c
T^{1/p},$$
hence the result follows. $\hfill{\Box}$
\end{section}


\begin{section}{Proof of the main theorem}\label{ch4:sec6}

At this stage, we are ready to present the proof of our main result
(Theorem \ref{ch4:mainresult}). This depends essentially on the passage
to the limit in the family of solutions $(\r^{\eps},\k^{\eps})$ of
system (\ref{ch4:app_model}). Since $\k^{\eps}_{x}\neq 0$, we multiply
the first equation of (\ref{ch4:app_model}) by $\k^{\eps}_{x}$ and we
rewrite system (\ref{ch4:app_model}) in
terms of $A^{\eps}$, we obtain:
\begin{equation}\label{ch4:sys_A}
\left\{
\begin{aligned}
& \k^{\eps}_{t}\k^{\eps}_{x}=\eps
\k^{\eps}_{x}\k^{\eps}_{xx}+\r^{\eps}_{x}A^{\eps}_{x}\quad &\mbox{on}&
\quad I_{T}\\
& \r^{\eps}_{t}=\eps \r^{\eps}_{xx}+A^{\eps}_{x}\quad &\mbox{on}&
\quad I_{T}.
\end{aligned}
\right.
\end{equation}
We will pass to the limit in the framework of viscosity solutions for
the first equation of (\ref{ch4:sys_A}), and in the distributional sense
for the second equation. We start with the following proposition.
\begin{proposition}\textit{\textbf{(Local uniform
      convergence)}}\label{ch4:unif_bound_A}\\
The sequences $(\r^{\eps})_{\eps}$, $(\r_{x}^{\eps})_{\eps}$,
$(\k^{\eps})_{\eps}$, $(A^{\eps})_{\eps}$ and
$(A^{\eps}_{x})_{\eps}$ converge (up to extraction of a subsequence)
locally uniformly in $I_{T}$ as $\eps$ goes to zero.
\end{proposition}
\noindent {\bf Proof.} Let $V$ be an open compactly contained subset of
$I_{T}$. The constants that will appear in the proof are all independent
of $\eps$. However, they may depend on other fixed parameters including
$V$. The idea is to give an $\eps$-uniform control of the modulus of continuity in space
and in time of the quantities mentioned in Proposition
(\ref{ch4:unif_bound_A}),
which gives the local uniform convergence.  The $\eps$-uniform control
on the space modulus of continuity will be derived from the Corollary
\ref{ch4:k_{x}} and Lemma \ref{ch4:control_x}, while the  $\eps$-uniform control
on the time modulus of continuity will be derived from Lemma
\ref{ch4:mod_con_t}. The proof is divided into five steps.\\

\noindent {\bf Step 1. }\textsf{(Convergence of $A^{\eps}$ and $A^{\eps}_{x}$)}\\

\noindent From (\ref{ch4:MT5}), we know that
$\left\|\frac{\r^{\eps}_{x}}{\k^{\eps}_{x}}\right\|_{\infty}\leq
1$. We apply the interior $L^{p}$, $p>1$, estimates for parabolic
equations (see for instance
\cite[Theorem 7.13, page 172]{lieberman}) to the term
$A^{\eps}$ satisfying (\ref{ch4:para_qA}), we obtain:
\begin{equation}\label{ch4:am3fhh}
\|A^{\eps}\|_{W^{2,1}_{p}(V)}\leq c_{2} \|A^{\eps}\|_{p,V'},
\end{equation}
where $V'$ is any open subset of $I_{T}$ satisfying $V\subset\!\subset
V'\subset\!\subset I_{T}$. The constant $c_{2}~=~c_{2}(p,\tau,V,V')$ can
be chosen independent of $\eps$ first by applying a parabolic rescaling
of (\ref{ch4:para_qA}), and then using the fact that the factor multiplied by
$A^{\eps}_{xx}$ in (\ref{ch4:para_qA}) satisfying $1\leq 1+\eps \leq
2$. At this point, we apply Lemma
\ref{ch4:3lcom:lem1} for $A^{\eps}$ on $V'$, we get:
\begin{equation}\label{ch4:am3fhh1}
\|A^{\eps}\|_{p,V'}\leq c_{3} (\|A^{\eps}\|_{1,I_{T}}+1),
\end{equation}
and hence the above two equations (\ref{ch4:am3fhh}) and
(\ref{ch4:am3fhh1}) give:
\begin{equation}\label{ch4:qA4}
\|A^{\eps}\|_{W^{2,1}_{p}(V)}\leq c_{4}(\|A^{\eps}\|_{1,I_{T}}+1).
\end{equation}
We estimate the right hand side of (\ref{ch4:qA4}) in the following way:
\begin{eqnarray*}
\|A^{\eps}\|_{1,I_{T}}&=&\int_{I_{T}}|\r^{\eps}_{x}-\tau\k^{\eps}|\\
&\leq& \int_{I_{T}}\k^{\eps}_{x}+\tau|\k^{\eps}|\\
&\leq& (2+\t)T,
\end{eqnarray*}
where we have used the fact that $|\r_{x}^{\eps}|<\k^{\eps}_{x}$ (see
(\ref{ch4:MT5}) of Theorem \ref{ch4:chap3mainresult}) in the second
line, and the fact that $0\leq |\k^{\eps}| \leq 1$ (see Remark
\ref{ch4:fshA}) in the last line. Therefore,
inequality (\ref{ch4:qA4}) implies:
\begin{equation}\label{ch4:qA5}
\|A^{\eps}\|_{W^{2,1}_{p}(V)}\leq c_{5},\quad 1<p<\infty.
\end{equation}
We use the above inequality for $p>3$. In this case, the Sobolev
embedding in Hölder spaces (see \cite[Lemma
2.8]{IJM_PI}) gives:
$$W^{2,1}_{p}(V)\hookrightarrow C^{1+\a,\frac{1+\a}{2}}(V),\quad \alpha=1-{3}/{p}$$
and hence (\ref{ch4:qA5}) implies:
\begin{equation}\label{ch4:prelim1}
\|A^{\eps}\|_{C^{1+\a,\frac{1+\a}{2}}(V)}\leq
c_{6},
\end{equation}
which guarantees the equicontinuity and the equiboundedness of
$(A^{\eps})_{\eps}$ and $(A^{\eps}_{x})_{\eps}$. By the
Arzela-Ascoli Theorem (see for instance \cite{brez1}), we finally obtain
\begin{equation}\label{ch4:lim1}
A^{\eps}\longrightarrow A \quad\mbox{and}\quad
A_{x}^{\eps}\longrightarrow A_{x},
\end{equation}
up to a subsequence, uniformly on $V$ as $\eps\rightarrow 0$.\\

\noindent {\bf Step 2. }\textsf{(Convergence of $\k^{\eps}$)}\\

\noindent We control the modulus of continuity of $\k^{\eps}$ in space
and in time, locally uniformly in $\eps$.\\

\noindent {\bf Step 2.1. }\textsf{(Control of the modulus of continuity
  in time)}\\

\noindent The first equation of (\ref{ch4:sys_A}) gives:
$$\k^{\eps}_{t}=\eps\k^{\eps}_{xx}+\frac{\r^{\eps}_{x}}{\k^{\eps}_{x}}A^{\eps}_{x},$$
and hence, using the fact that $\left\|\frac{\r^{\eps}_{x}}{\k^{\eps}_{x}}\right\|_{\infty}\leq
1$, together with (\ref{ch4:prelim1}), we get:
\begin{equation}\label{ch4:3b3al}
\|\k^{\eps}_{t}-\eps\k^{\eps}_{xx}\|_{\infty,V}\leq
\left\|\frac{\r^{\eps}_{x}}{\k^{\eps}_{x}}\right\|_{\infty,V}\|A_{x}\|_{\infty,V}\leq c_{6}.
\end{equation}
Also, by (\ref{ch4:fshA1}), we have:
$$\|\k^{\eps}\|_{\infty,V}\leq 1.$$
This uniform bound on $\k^{\eps}$ together with (\ref{ch4:3b3al}) permit
to use Lemma \ref{ch4:mod_con_t} to conclude that
\begin{equation}\label{ch4:lim2}
|\k^{\eps}(x,t)-\k^{\eps}(x,t+h)|\leq c_{7}h^{\beta} ,
\quad(x,t), (x,t+h)\in V,\quad 0<\beta<1,
\end{equation}
which controls the modulus of continuity of
$\k^{\eps}$ with respect to $t$ uniformly in $\eps$. We now move to
control the modulus of continuity in space.\\

\noindent {\bf Step 2.2 }\textsf{(An $\eps$-uniform bound on $S(0)$)}\\

\noindent Recall the definition (\ref{ch4:S(t)}) of $S(t)$:
$$S(t)=\int_{I}\sum_{\pm}\theta^{\pm,\eps}(x,t)\log
\theta^{\pm,\eps}(x,t) dx,$$
with
$$\theta^{\pm,\eps}=\frac{\k^{\eps}_{x}\pm \r^{\eps}_{x}}{2}.$$
Hence
$$S(0)=\int_{I}\frac{\k^{0,\eps}_{x}+
  \r^{0,\eps}_{x}}{2}\log\left(\frac{\k^{0,\eps}_{x}+ \r^{0,\eps}_{x}}{2}\right)
+\int_{I}\frac{\k^{0,\eps}_{x}-\r^{0,\eps}_{x}}{2}
\log\left(\frac{\k^{0,\eps}_{x}-\r^{0,\eps}_{x}}{2}\right).$$
Using the elementary identities $x \log x \leq x^{2}$ and $(x\pm
y)^{2}\leq 2(x^{2}+y^{2})$, we compute:
\begin{eqnarray}\label{ch4:shubdi3lik}
S(0)&\leq& \int_{I}\left(\frac{\k^{0,\eps}_{x}+
    \r^{0,\eps}_{x}}{2}\right)^{2}+\int_{I}\left(\frac{\k^{0,\eps}_{x}-
    \r^{0,\eps}_{x}}{2}\right)^{2}\nonumber\\
&\leq& \|\r^{0,\eps}_{x}\|^{2}_{2,I}+\|\k^{0,\eps}_{x}\|^{2}_{2,I}.
\end{eqnarray}
From (\ref{ch4:cond_init_ap2}), we know
that:
$$
|\r^{0,\eps}_{x}|=\left|\frac{\r^{0}_{x}+\eps \tau
    \phi^{'}}{(1+\eps)^{2}}\right|\leq
\frac{|\r^{0}_{x}|+\eps}{(1+\eps)^{2}}\leq |\r^{0}_{x}|+1,
$$
and
$$
|\k^{0,\eps}_{x}|=\left|\frac{\k^{0}_{x}+\eps}{1+\eps}\right|\leq |\k^{0}_{x}|+1.
$$
Using the above two inequalities into (\ref{ch4:shubdi3lik}), we deduce
that:
$$S(0)\leq 2(\|\r^{0}_{x}\|^{2}_{2,I}+\|\k^{0}_{x}\|^{2}_{2,I}+2).$$

\noindent {\bf Step 2.3. }\textsf{(Control of the modulus of continuity
  in space and conclusion)}\\

\noindent We use the uniform bound obtained for $S(0)$ in Step 2.1,
together with the special control (\ref{ch4:cont_k_x}) of
$\k^{\eps}_{x}$ given in Corollary \ref{ch4:k_{x}}, we get for all
$0\leq t \leq T$:
$$
\int_{I}\k^{\eps}_{x}(x,t)\log(\k^{\eps}_{x}(x,t))dx\leq
2(\|\r^{0}_{x}\|^{2}_{2,I}+\|\k^{0}_{x}\|^{2}_{2,I}+2)+\frac{\tau^{2}T}{2}+2,
$$
therefore
\begin{equation}\label{ch4:hasaTo}
\int_{I}\k^{\eps}_{x}(x,t)\log(\k^{\eps}_{x}(x,t))dx\leq c_{8},\quad
\forall\,0\leq t\leq T.
\end{equation}
Inequality (\ref{ch4:hasaTo}) permits to use Lemma \ref{ch4:control_x},
hence we obtain:
\begin{equation}\label{ch4:lim3}
|\k^{\eps}(x+h,t)-\k^{\eps}(x,t)|\leq \frac{c_{9}}{|\log h|},\quad (x,t), (x+h,t)\in
I_{T},
\end{equation}
Inequalities
(\ref{ch4:lim2}) and (\ref{ch4:lim3}) give the equicontinuity of the
sequence $(\k^{\eps})_{\eps}$ on $V$, and again by the Arzela-Ascoli
Theorem, we get:
\begin{equation}\label{ch4:lim4}
\k^{\eps}\rightarrow \k,
\end{equation}
up to a subsequence, uniformly on $V$ as $\eps\rightarrow 0$.\\

\noindent {\bf Step 3. }\textsf{(Convergence of $\r^{\eps}$)}\\

\noindent As in step 2, we control the modulus of continuity of
$\r^{\eps}$ in space and in time, locally uniformly with respect to
$\eps$.\\

\noindent {\bf Step 3.1. }\textsf{(Control of the modulus of continuity
  in time)}\\

\noindent The second equation of (\ref{ch4:sys_A}) gives:
$$\r^{\eps}_{t}-\eps \r^{\eps}_{xx}=A_{x}^{\eps},$$
hence, from (\ref{ch4:prelim1}), we deduce that:
$$\|\r^{\eps}_{t}-\eps \r^{\eps}_{xx}\|_{\infty,V}\leq c_{6},$$
and from (\ref{ch4:fshA1}), we have:
$$\|\r^{\eps}\|_{\infty,V}\leq 1.$$
The above two inequalities permit to use Lemma \ref{ch4:mod_con_t}, we
finally get:
\begin{equation}\label{ch4:lim2forr}
|\r^{\eps}(x,t)-\r^{\eps}(x,t+h)|\leq c_{7}h^{\beta} ,
\quad(x,t), (x,t+h)\in V,\quad 0<\beta<1,
\end{equation}
which controls the modulus of continuity of
$\r^{\eps}$ with respect to $t$ uniformly in $\eps$.\\

\noindent {\bf Step 3.2. }\textsf{(Control of the modulus of continuity
  in space and conclusion)}\\

\noindent The control of the space modulus of continuity is based on the
following observation. From (\ref{ch4:MT5}), we know that
$|\r^{\eps}_{x}|\leq \k^{\eps}_{x}$ on $I_{T}$. Using this inequality, we
get, for every $(x,t),
(x+h,t)\in I_{T}$:
$$|\r^{\eps}(x+h,t)-\r^{\eps}(x,t)|\leq \int_{x}^{x+h}|\r^{\eps}_{x}(y,t)|dy
\leq \int_{x}^{x+h}\k^{\eps}_{x}(y,t)dy\leq
|\k^{\eps}(x+h,t)-\k^{\eps}(x,t)|.$$
Inequality (\ref{ch4:lim3}) gives immediately that:
\begin{equation}\label{ch4:limlam}
|\r^{\eps}(x+h,t)-\r^{\eps}(x,t)|\leq \frac{c_{9}}{|\log h|},\quad (x,t), (x+h,t)\in
I_{T}.
\end{equation}
From (\ref{ch4:lim2forr}) and (\ref{ch4:limlam}), we deduce that:
\begin{equation}\label{ch4:lim5}
\r^{\eps}\rightarrow \r,
\end{equation}
up to a subsequence, uniformly on $V$ as $\eps\rightarrow 0$.\\

\noindent {\bf Step 4. }\textsf{(Convergence of $\r_{x}^{\eps}$ and conclusion)}\\

\noindent In fact, this follows from Step 1, Step 2, and
the fact that
\begin{equation}\label{ch4:lim6}
\r_{x}^{\eps}=A^{\eps}+\t\k^{\eps}\rightarrow \r_{x},
\end{equation}
uniformly on $V$ as $\eps\rightarrow 0$. In this case, we also deduce that
$$
A=\r_{x}-\tau\kappa.
$$
The proof of Proposition \ref{ch4:unif_bound_A} is done.$\hfill{\Box}$\\

We now move to the proof of the main result.\\

\noindent {\bf Proof of Theorem \ref{ch4:mainresult}.} We first remark
that $\k^{\eps}$ is a viscosity solution of the first equation of (\ref{ch4:sys_A}):
\begin{equation}\label{ch4:am}
\k^{\eps}_{t}\k^{\eps}_{x}-\eps\k^{\eps}_{x}\k^{\eps}_{xx}-\r^{\eps}_{x}A^{\eps}_{x}=0\quad
\mbox{on}\quad I_{T}.
\end{equation}

\begin{rem}\label{ch4:rem4}
The equation (\ref{ch4:am}) can be viewed as the following
Hamilton-Jacobi equation of second order:
\begin{equation}\label{ch4:Hamil1}
H^{\eps}(X,D\k^{\eps},D^{2}\k^{\eps})=0,\quad X=(x,t)\in I_{T}
\end{equation}
with
$$D\k^{\eps}=(\k^{\eps}_{x},\k^{\eps}_{t})\quad\mbox{and}\quad D^{2}\k^{\eps}=\left(
\begin{array}{cc}
\k^{\eps}_{xx} & \k^{\eps}_{xt} \\
\k^{\eps}_{tx} & \k^{\eps}_{tt} \\
\end{array}
\right),
$$
where $H^{\eps}$ is the Hamiltonian function given by:
\begin{equation}\label{ch4:Hamil2}
\begin{aligned}
H^{\eps}: \;&I_{T}\times\R^{2}\times
M^{2\times2}_{sym}\hspace{-0.2cm}&
&\longrightarrow&&\hspace{-0.1cm}\R&\\
&\hspace{0.5cm}(X,p,M)&&\longmapsto&
&\hspace{-0.1cm}H^{\eps}(X,p,M)=p_{1}p_{2}-\eps
p_{1}M_{11}-\r^{\eps}_{x}(X)A^{\eps}_{x}(X),&
\end{aligned}
\end{equation}
$p=(p_{1},p_{2})$ and $M=(M_{ij})_{i,j=1,2}.$
\end{rem}
From (\ref{ch4:lim1}) and (\ref{ch4:lim6}), we deduce that
$(H^{\eps})_{\eps}$ converges locally uniformly in
$I_T\times\R^{2}\times M^{2\times 2}_{sym}$ to the function $H$
given by:
\begin{equation}\label{ch4:Hamil_f}
\begin{aligned}
H: \;&I_{T}\times\R^{2}\times
M^{2\times2}_{sym}\hspace{-0.2cm}&
&\longrightarrow&&\hspace{-0.1cm}\R&\\
&\hspace{0.5cm}(X,p,M)&&\longmapsto&
&\hspace{-0.1cm}H(X,p,M)=p_{1}p_{2}-\r_{x}(X)A_{x}(X).&
\end{aligned}
\end{equation}
This, together with the local uniform convergence of $\k^{\eps}$ to $\k$ (see
\ref{ch4:lim4}), and the fact that $\k^{\eps}$ is a
viscosity solution of (\ref{ch4:am}), permit to use the stability of
viscosity solutions (see Theorem
\ref{ch4:stability}), which proves that $\k$ is a viscosity solution of
\begin{equation}\label{ch4:9}
H(X,D\k,D^{2}\k)=\k_{t}\k_{x}-\r_{x}A_{x}=0\quad\mbox{in}\quad
I_{T}.
\end{equation}
We now pass to the limit $\eps\rightarrow 0$ in the second equation
of (\ref{ch4:sys_A}), we obtain
\begin{equation}\label{ch4:10}
\r_{t}=A_{x}\in C(I_{T})\quad\mbox{in}\quad \mathcal{D}'(I_{T}).
\end{equation}
From (\ref{ch4:9}) and (\ref{ch4:10}), we get:
\begin{enumerate}
\item $\k$ is a viscosity solution of $\k_{t}\k_{x}=\r_{t}\r_{x}$ in
$I_{T}$;
\item $\r$ is a distributional solution of
$\r_{t}=\r_{xx}-\tau\k_{x}$ in $I_{T}$.
\end{enumerate}
Inequality (\ref{ch4:7}) can be easily obtained by testing against
a nonnegative function $\phi\in C^{\infty}_{0}(I_T)$. Finally, let
us show how to retrieve the initial and boundary conditions.
Indeed, the local uniform convergence
$(\r^{\eps},\k^{\eps})\rightarrow (\r,\k)$, together with the
uniform control of the modulus of continuity of these solutions:
\begin{itemize}
  \item with respect to $x$ near $\partial I\times [0,T]$ by (\ref{ch4:lim3});
  \item with respect to $t$ near $I\times\{t=0\}$, away from $0$ and
  $1$ by (\ref{ch4:lim2}),
\end{itemize}
and the fact that $\k^{0,\eps}\rightarrow \k^{0}$,
$\r^{0,\eps}\rightarrow \r^{0}$ uniformly in $\bar
I$, show that $(\r,\k)\in (C(\overline{I_{T}}))^{2}$, so the initial
and boundary conditions are satisfied pointwisely, and the proof
of the main result is done. $\hfill{\Box}$
\end{section}




\begin{section}{Application: simulations for the evolution of
    elastoviscoplatic materials}\label{ch4:Simulation}
Motivated by the simulation of the elastoviscoplastic behavior
that are formulated by the model of Groma, Csikor and Zaiser
\cite{GCZ}, this section is devoted to write down the equations of
the displacement vector $u$ inside the crystal when a constant
exterior shear stress $\t$ (in this section we take $\t>0$) is applied
on the boundary walls (see
Figure \ref{IMP_fig}). Also, at the end of this section, we
present some numerical simulations revealing the evolution of a
crystal. Here, as we have already mentioned in the introduction, we suppose
that the distribution of dislocations is invariant by translation
in the $y$-direction. We consider a 2-dimensional crystal (Figure
\ref{IMP_fig}) with the displacement vector:
$$u=(u_{1},u_{2}):\R^{2}\longmapsto \R^{2}.$$
We rename the coordinates $(x,y)$ as $(x_{1},x_{2})$, and we call
$(e_{1},e_{2})$ the corresponding orthonormal basis. We
define the total strain by:
$$\eps_{ij}(u)=\frac{1}{2}\left(\partial_{j}u_{i}+\partial_{i}u_{j}\right)
\quad\mbox{with}\quad\partial_{j}u_{i}=\frac{\partial u_{i}}{\partial
x_{j}},\quad i,j=1,2.$$ This total strain could be decomposed into
two parts as follows:
\begin{equation}\label{Heur:lala7z1}
\eps_{ij}(u)=\eps^{e}_{ij}(u)+\eps_{ij}^{p},
\end{equation}
where $\eps^{e}(u)$ is the elastic strain and $\eps^{p}$ is the
plastic strain which is given by:
\begin{equation}\label{Heur:plas_str}
\eps_{ij}^{p}=\rho\eps_{ij}^{0},
\end{equation}
with
$$\eps_{ij}^{0}=\frac{1}{2}(1-\delta_{ij}),$$
in the special case of a single slip system where dislocations
move following the Burgers vector $\vec{b}=e_{1}$ ($\delta_{ij}$
is the usual Kronecker symbol). Here $\rho$ is the resolved
plastic strain. The stress field $\sigma$ inside the crystal is
given by:
\begin{equation}\label{Heur:PUit}
\sigma_{ij}=2\mu \eps_{ij}^{e}(u)+\lambda \delta_{ij}\left(
\sum_{k=1,2}\eps_{kk}^{e}(u)\right),
\end{equation}
with $\lambda, \mu >0$ are the constants of Lamé coefficients of
the crystal that are assumed (for simplification) to be isotropic.
This stress field $\sigma$ has to satisfy the equation of
elasticity:
\begin{equation}\label{eqn_of_elasticity}
\sum_{j=1,2} \frac{\partial \sigma_{ij}}{\partial x_j} = 0.
\end{equation}
Finally, the functions $\r,\k$ (solutions of (\ref{ch4:sys_rho_kappa}))
and $u$ are solutions of the following coupled system:
\begin{equation}\label{The_Full_SYS}
\left\{
\begin{aligned}
&\sum_{j=1,2} \frac{\partial \sigma_{ij}}{\partial x_j}
 &=\quad& 0 \quad &\mbox{on}&\quad I\times (0,\infty),\\
&\sigma_{ij} &=\quad& 2\mu \eps_{ij}^{e}(u)+\lambda
\delta_{ij}\left( \sum_{k=1,2}\eps_{kk}^{e}(u)\right)
&\mbox{on}&\quad I\times
(0,\infty),\\
&\eps_{ij}^{e}&=\quad&
\frac12\left(\partial_{j}u_{i}+\partial_{i}u_{j}\right)-\rho\eps^{0}_{ij}\quad
&\mbox{on}&\quad
I\times (0,\infty),\\
& \k_{t}\k_{x} &=\quad& \r_{t}\r_{x}\quad &\mbox{on}&\quad
I\times (0,\infty),\\
&\r_{t}&=\quad& \r_{xx}-\tau \k_{x}\quad &\mbox{on}&\quad
I\times (0,\infty),
\end{aligned}
\right.
\end{equation}
Equation (\ref{eqn_of_elasticity}) can be reformulated as:
\begin{equation}\label{Heur:UinsideI}
\left\{\begin{array}{ll}(\lambda+2\mu) \Delta
u_1+(\lambda+\mu)\partial_{12}u_2=0\\ \\
(\lambda+2\mu) \Delta
u_2+(\lambda+\mu)\partial_{21}u_1=\mu\partial_1\rho.
\end{array}
\right.
\end{equation}
Here $\partial_{2}\r=0$ is due to the homogeneity of the distribution of
dislocations in the $e_{2}$-direction.\\

\noindent \textbf{Calculation of $u$.} We first calculate the
value of the displacement $u$ on the boundary walls. Remark first
that since we are applying a constant shear stress field on the
walls, the stress field $\sigma$ there can be evaluated as:
$\sigma\cdot n=\pm \tau e_{2}$, $n=\pm e_{1}$, for $x=\pm1$,
\begin{equation}\label{Heur:siGbon}
\sigma=\left(
\begin{array}{cc}
0 & \tau\\
\tau & 0
\end{array}
\right),\quad \mbox{on}\quad \partial I.
\end{equation}
Using (\ref{Heur:siGbon}) and (\ref{Heur:PUit}), we can derive the
following equations on the boundary:
\begin{equation}\label{Heur:eqnn1}
\left\{
\begin{aligned}
&\partial_{1}u_{1}=0 \quad &\mbox{on}& \quad \partial I,\\
&\mu (\partial_{1}u_{2}-\r)=\tau \quad &\mbox{on}& \quad \partial I.
\end{aligned}
\right.
\end{equation}
Equation (\ref{Heur:UinsideI}) leads to the following two equations
inside $I$:
\begin{equation}\label{Heur:eqnn2}
\left\{
\begin{aligned}
& \partial_{1}[(\lambda + 2 \mu) \partial_{1}u_{1}]=0 \quad &\mbox{on}&
\quad I\\
&\partial_{1}(\partial_{1}u_{2}-\r)=0\quad &\mbox{on}& \quad I.
\end{aligned}
\right.
\end{equation}
Combining (\ref{Heur:eqnn1}) and (\ref{Heur:eqnn2}), with
$$u_{1}(0,x_{2})=u_{2}(0,x_{2})=0,$$
finally lead to:
\begin{equation}\label{Heur:eqnn4}
\left\{
\begin{aligned}
& u_{1}(x_{1},x_{2})=0, \quad &(x_{1},x_{2})\in I\times \R&\\
& u_{2}(x_{1},x_{2})=\frac{\tau}{\mu} x_{1} +
\int^{x_{1}}_{0}\r(x)dx,\quad &(x_{1},x_{2})\in I\times \R.&
\end{aligned}
\right.
\end{equation}
\noindent \textbf{Formal computation of the long time solution of system
  (\ref{ch4:sys_rho_kappa}).} All computations that will be done here
are formal. We seek to calculate long time solutions for system
(\ref{ch4:sys_rho_kappa}). For this reason, we first calculate the long
time solution (stationary solution, $\k^{\eps}_{t}=\r_{t}^{\eps}=0$) of
the regularized $\eps$-system (\ref{ch4:app_model}), and then
we pass to the limit $\eps\rightarrow 0$. Doing some computations, the
long time $(t\rightarrow \infty)$ solution for the $\eps$-system is given by:
\begin{equation}\label{l_t_c_f}
\r^{\eps}(x,t)= B \left(\cosh\left(\frac{\tau x}{1+\eps} \right)-
  \cosh\left(\frac{\tau}{1+\eps}\right) \right)\quad \mbox{and}\quad
\k^{\eps}(x,t)=B \sinh \left(\frac{\tau}{1+\eps}\right),
\end{equation}
with $B=1/\sinh\left(\frac{\tau}{1+\eps}\right)$. Passing (again
formally) to the limit as $\eps \rightarrow 0$ in (\ref{l_t_c_f}), and
using (\ref{Heur:eqnn4}), we can compute the long
time displacement $u_{2}$ inside the material. In fact, we have:
\begin{equation}\label{l_t_c_f_f_f}
u_{2}(x_{1},x_{2})=\left(\frac{\tau}{\mu}-\frac{\cosh \tau}{\sinh
    \tau}\right)x_{1} + \frac{\sinh \tau x_{1}}{\tau \sinh \tau}.
\end{equation}

\noindent \textbf{Numerical simulations.} The displacement $u$ is
numerically computed by discretizing system (\ref{ch4:sys_rho_kappa})
using an upwind numerical scheme. The space and time steps in this
scheme are well chosen in order to satisfy a CFL condition (for the
details, see \cite[Appendix]{hassan_these}). In Figure \ref{p:hif5}, we
show successively the initial state of the crystal
at time $t=0$ without any applied stress, then the instantaneous (elastic) deformation of
the crystal when we apply the shear stress $\tau>0$ at time $t=0^+$. The
deformation of the crystal evolves in time and finally converges
numerically to some particular deformation which is shown on the last
picture after a very long time. This kind of behavior is called
elasto-visco-plasticity in mechanics.
Moreover, on the last picture, we observe the presence of boundary layer
deformations. This effect is directly related to the introduction of the
back stress $\tau_b=\frac{\th^{+}_{x}-\th^{-}_{x}}{\th^{+}+\th^{-}}$ in the model.

\begin{figure}[!h]
\begin{center}
\begin{tabular}{lcr}
\epsfig{figure=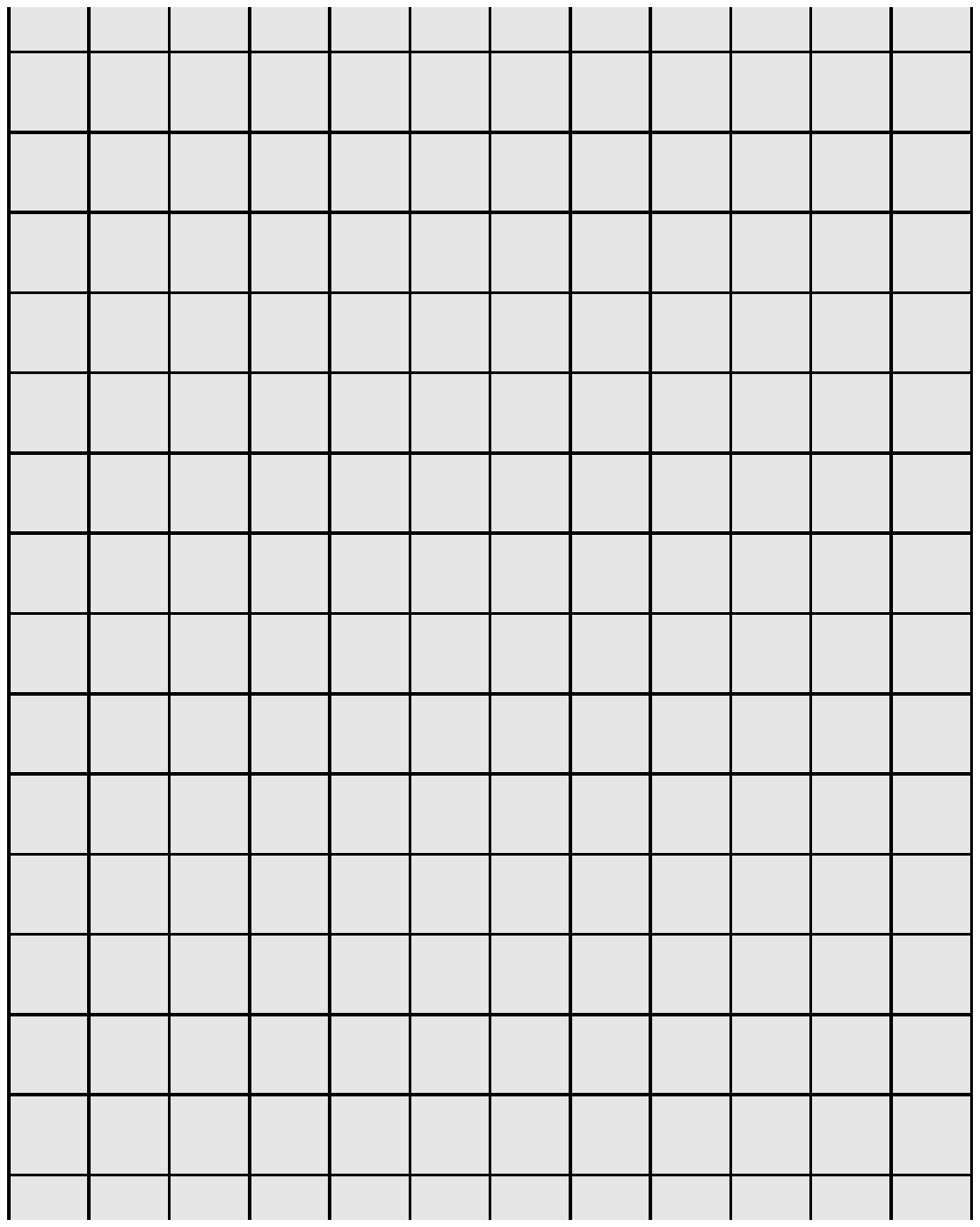,width=25.5mm}
&$\quad \quad$
\epsfig{figure=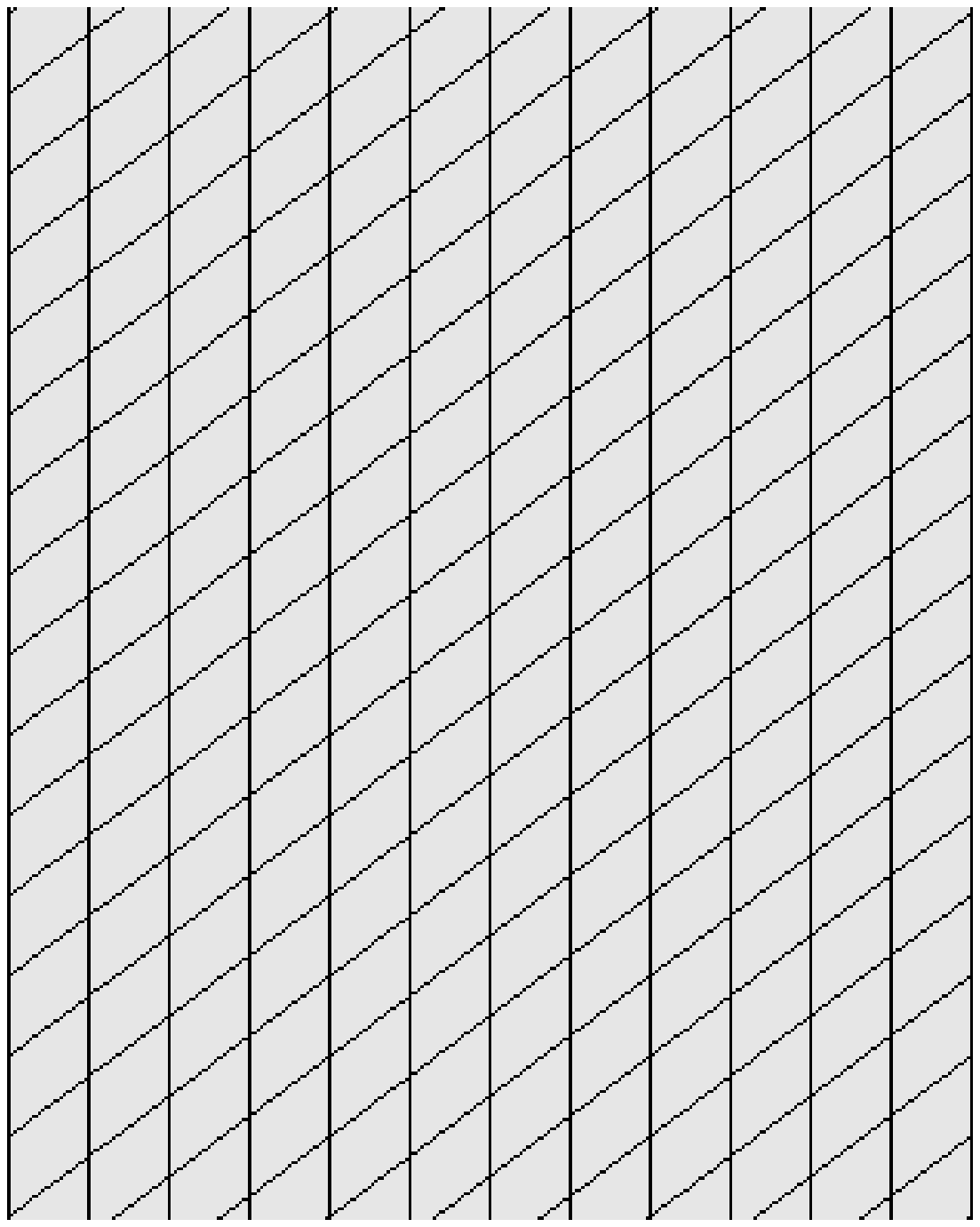,width=25.5mm}
&$\quad \quad$
\epsfig{figure=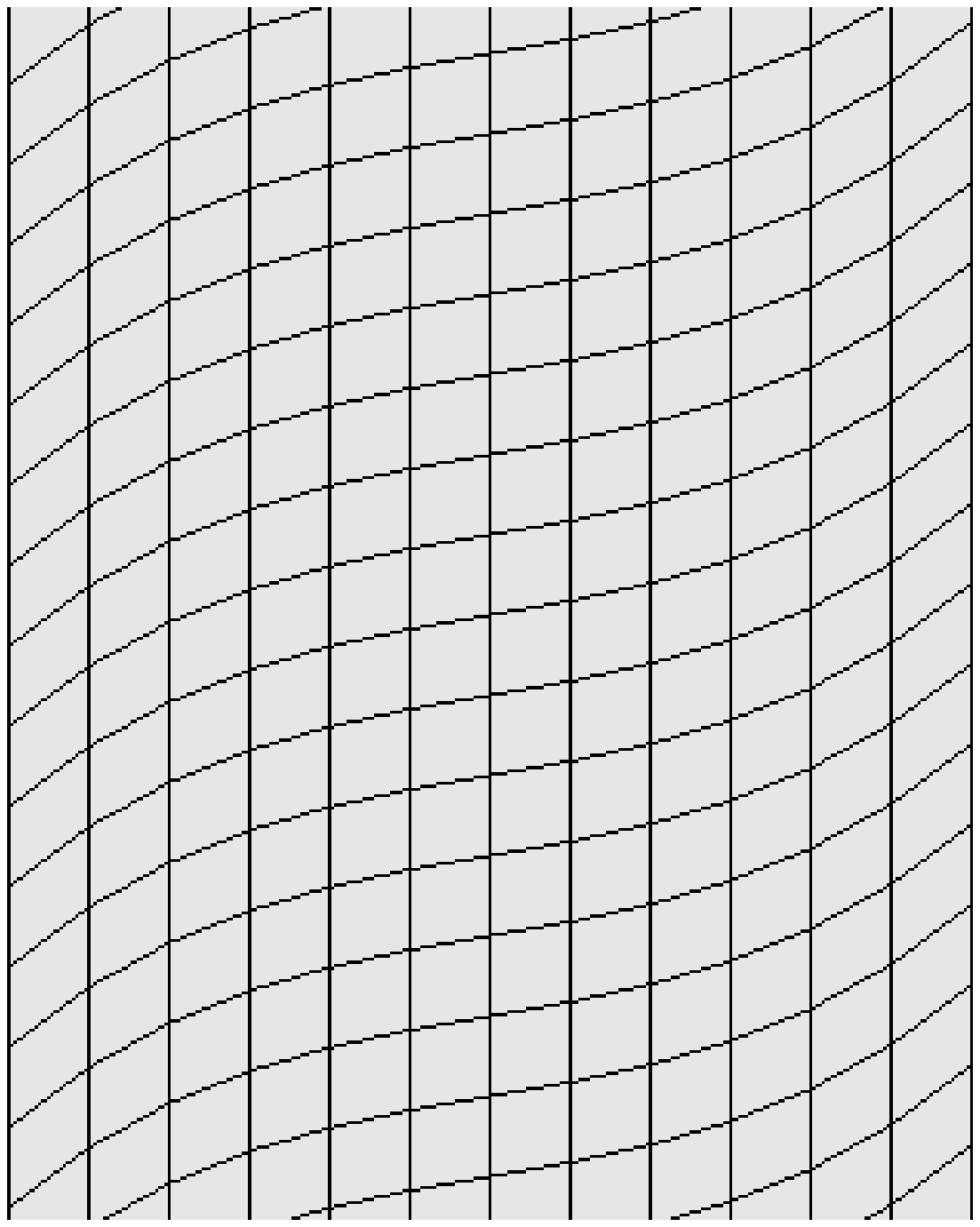,width=25.5mm}\\
$\hspace{0.5cm}$ a) $t=0^{-}$ & $\quad\quad$ b) $t=0^{+}$ & c) $t=+\infty$
$\hspace{0.33cm}$
\end{tabular}
\end{center}
\caption{Deformation of a slab for model (\ref{ch4:sys_rho_kappa}).}\label{p:hif5}
\end{figure}
\begin{rem}
In the last image of Figure \ref{p:hif5}, it seems that the long time
displacement $u_{2}$ inside the material is coherent with our formal
computation for the term $u_{2}$ (see equation (\ref{l_t_c_f_f_f})).
\end{rem}

\end{section}


\begin{section}{Appendix}\label{ch4:secA}

\noindent {\bf A1. Proof of Lemma \ref{ch4:mod_con_t} (control of the
  modulus of continuity in time)}\\

\noindent Let $V$ be a compactly contained subset of $I_{T}$. Throughout the
proof, the constant $c$ may take several values but only depending
on $V$. Since $V\subset\!\subset I_{T}$, then there is a rectangular
cube of the form
$$\mathcal{Q}=(x_{1},x_{2})\times (t_{1},t_{2}),$$
such that $V \subset\!\subset \mathcal{Q} \subset\!\subset I_{T}$.
In this case, there exists a constant $\eps_{0}$, also depending on
$V$ such that for any $$0<\eps<\eps_{0},$$ and any $(x,t)\in V$, we
have:
$$
(x-2\sqrt{\eps},x+2\sqrt{\eps})\times\{t\}\subset \mathcal{Q}.
$$
Moreover, for any $(x,t),(x,t+h)\in V$, we can always find two
intervals $\mathcal{I}$ and $\mathcal{J}$ such that
$$(t,t+h)\subset \mathcal{I}\subset\!\subset \mathcal{J}, \mbox{
with } \{x\}\times \mathcal{I}\subset
\mathcal{Q}\quad\mbox{and}\quad\{x\}\times \mathcal{J}\subset
\mathcal{Q}.$$ Let us indicate that these intervals might have
different lengths depending on $h$ and $V$ but we always have
$$
|\mathcal{J}|, |\mathcal{I}|\leq |t_{2}-t_{1}|.
$$
Consider the following rescaling of the function $u^\eps$ defined
by:
\begin{equation}\label{ch4:Vd0}
\tilde{u}^{\eps}(x,t)=u^{\eps}(\sqrt{\eps}x,t).
\end{equation}
This function satisfies
$$\tilde{u}^{\eps}_{t}=\tilde{u}^{\eps}_{xx}+\tilde{f}^{\eps},\quad (x,t)\in
(0,1/\sqrt{\eps})\times (0,T),$$ where
$\tilde{f}^{\eps}(x,t)=f^{\eps}(\sqrt{\eps}x,t)$. Take
$(x_{0},t_{0})$, $(x_{0},t_{0}+h)$ in $V$, and let
$$\mathcal{Q}_{1}=(x_{0}-\sqrt{\eps},x_{0}+\sqrt{\eps})\times
\mathcal{I}\quad\mbox{and}\quad
\mathcal{Q}_{2}=(x_{0}-2\sqrt{\eps},x_{0}+2\sqrt{\eps})\times
\mathcal{J}.$$ These two cylinders are transformed by the above
rescaling into
$$\tilde{\mathcal{Q}}_{1}=\left(\frac{x_{0}}{\sqrt{\eps}}-1,\frac{x_{0}}{\sqrt{\eps}}+1\right)
\times \mathcal{I}\quad\mbox{and}\quad
\tilde{\mathcal{Q}}_{2}=\left(\frac{x_{0}}{\sqrt{\eps}}-2,\frac{x_{0}}{\sqrt{\eps}}+2\right)\times
\mathcal{J}.$$ We apply the interior $L^{p}$, $p>3$, estimates for parabolic
equations (see for instance
\cite[Theorem 7.13, page 172]{lieberman}) to the function $\tilde{u}^{\eps}$ over
the domains $\tilde{\mathcal{Q}}_{1}\subset\!\subset
\tilde{\mathcal{Q}}_{2}$, we get
\begin{equation}\label{ch4:Vd3}
\|\tilde{u}^{\eps}\|_{W^{2,1}_{p}(\tilde{\mathcal{Q}}_{1})}\leq c
(\|\tilde{u}^{\eps}\|_{p,\tilde{\mathcal{Q}}_{2}}+\|\tilde{f}^{\eps}\|
_{p,\tilde{\mathcal{Q}}_{2}}).
\end{equation}
Using the local $\eps$-uniform boundedness of $(u^{\eps})_{\eps}$, and
$(f^{\eps})_{\eps}$, we get:
$$
\|\tilde{u}^{\eps}\|^{p}_{L^{p}(\tilde{\mathcal{Q}}_{2})}\leq c\quad
\mbox{and}\quad \|\tilde{f}^{\eps}\|^{p}_{L^{p}(\tilde{\mathcal{Q}}_{2})}\leq c.
$$
hence, inequality (\ref{ch4:Vd3}) implies:
\begin{equation}\label{ch4:Vd6}
\|\tilde{u}^{\eps}\|_{W^{2,1}_{p}(\tilde{\mathcal{Q}}_{1})}\leq c.
\end{equation}
We use the Sobolev embedding in Hölder spaces (see for instance
\cite[Lemma 2.8]{IJM_PI}):
$$
W^{2,1}_{p}(\tilde{\mathcal{Q}}_{1})\hookrightarrow
C^{{1+\alpha}{\frac{1+\alpha}{2}}}(\tilde{\mathcal{Q}}_{1}),\quad p>3,\; \alpha=1-{3}/{p},
$$
to obtain, from (\ref{ch4:Vd6}), that:
$$\|\tilde{u}^{\eps}\|_{C^{{1+\alpha}{\frac{1+\alpha}{2}}}(\tilde{\mathcal{Q}}_{1})}\leq
c,$$ and hence
$$\frac{|\tilde{u}^{\eps}(x_{0}/\sqrt{\eps},t_{0}+h)-\tilde{u}^{\eps}(x_{0}/\sqrt{\eps},t_{0})|}
{h^{\frac{1+\alpha}{2}}}\leq c,$$ then from (\ref{ch4:Vd0}),
$$\frac{|{u}^{\eps}(x_{0},t_{0}+h)-{u}^{\eps}(x_{0},t_{0})|}{h^{\frac{1+\alpha}{2}}}\leq
c.$$ Choosing $\beta=\frac{1+\alpha}{2}$ we get the desired
 result.$\hfill{\Box}$\\

\noindent {\bf A2. Proof of Lemma \ref{ch4:vnlemma} (An interior estimate
  for the heat equation)}\\

\noindent Recall that $a$ is a solution of the heat equation on $I_{T}$,
$$a_{t}=a_{xx}.$$
The proof of lemma \ref{ch4:vnlemma} is a direct computation using
a mean value formula for solutions of the heat equations. Usually,
basic mean value formulae of the solution of the heat equation are
expressed through unbounded kernels (see for example \cite[Theorem
1]{FaGa}), where $a$ can be expressed as:
\begin{equation}\label{ch4:mean_basic}
a(x_{0},t_{0})=(4\pi r^{2})^{-1/2}\int_{\Omega_{r}(x_{0},t_{0})} a(x,t) \frac{(x_{0}-x)^{2}}{4(t_{0}-t)^{2}}dxdt.
\end{equation}
Here, $(x_{0},t_{0})\in I_{T}$, $(x,t)\in
\Omega_{r}(X_{0})$, and $r>0$ small enough in order to
ensure that the parabolic ball of radius $r$:
\begin{equation}\label{ch4:def_O}
\Omega_{r}(x_{0},t_{0})=\left\{(x,t);\; t_{0}-r^{2}<t<t_{0},\;\;
  (x-x_{0})^{2}< 2(t_{0}-t)\log\left(\frac{r^{2}}{t_{0}-t}\right)
\right\}\subset I_{T}.
\end{equation}
In our case, we need a mean value formula similar to
(\ref{ch4:mean_basic}) but with a bounded kernel on
$\Omega_{r}(x_{0},t_{0})$. In \cite{Garo}, the authors have given
such a representation formula for the solution of the heat
equation. The following is a direct corollary of \cite[Theorem
3.1]{Garo}):

\begin{corollary}\textit{\textbf{(Mean value formula with bounded kernels,
      \cite[Theorem 3.1]{Garo})}}\label{ch4:mean_value_cor}\\
Let $u\in C^{2}(\mathcal{D})$ be a solution of the heat equation:
$$u_{t}=u_{xx}\quad \mbox{on} \quad \mathcal{D},$$
where $\mathcal{D}$ is an open subset of $\R^{2}$ containing the
modified unit parabolic ball $\Omega'_{r}(x_{0},t_{0})$, $r>0$,
with
$$\Omega'_{r}(x_{0},t_{0})=\left\{(x,t);\; t_{0}-r^{2}<t<t_{0},\;\;
  |x-x_{0}|^{2}< 8(t_{0}-t)\log \left(\frac{r^{2}}{t_{0}-t}\right)
\right\}.$$ Then we have:
\begin{equation}\label{ch4:G_exp_Bnd_any}
u(x_{0},t_{0})=\frac{\bar{c}}{|\Omega'_{r}(x_{0},t_{0})|}\int_{\Omega'_{r}(x_{0},t_{0})}
u(x,t)E\left(\frac{x-x_{0}}{r},\frac{t-t_{0}}{r^{2}}\right)dxdt,
\end{equation}
where $\bar{c}>0$, $|\Omega'_{r}(x_{0},t_{0})|=\bar{c}r^{3}$, and
the kernel $E$ satisfies:
\begin{equation}\label{ch4:G_exp_Bnd}
\|E(x,t)\|_{\infty, \Omega'_{1}(0,0)}\leq c,
\end{equation}
and $c>0$ is a fixed positive constant.
\end{corollary}

\begin{rem}
In the above corollary, which is an application of \cite[Theorem
3.1]{Garo} in the case $m=3$, an explicit expression of $E$ can be
given by:
$$E(x,t)=\frac{\omega_{3}}{16
  \pi^{2}}\left(-x^{2}+8t\log(-t)\right)^{3/2}\left[\frac{x^{2}}{4t^{2}}+
\frac{3(-x^{2}+8t\log(-t))}{20t^{2}}\right],$$
where $\omega_{3}$ is the volume of the unit ball in $\R^{3}$. For a
more general expression of $E$, we send the reader to \cite[Equality
(3.6) of Theorem 3.1]{Garo}.
\end{rem}

\end{section}

\noindent \textbf{Acknowledgments}. This work was supported by the
contract ANR MICA (2006-2009). The authors would like to thank
J\'{e}r\^{o}me Droniou for his valuable remarks while reading the
manuscript of the paper.

\bibliographystyle{siam}
\bibliography{biblio}

\end{document}